\PassOptionsToPackage{hyphens}{url}
\documentclass[11pt]{article}
\usepackage[margin=1in]{geometry}
\usepackage[T1]{fontenc}
\usepackage[utf8]{inputenc}
\usepackage{lmodern}
\usepackage{microtype}
\usepackage{amsmath,amssymb,amsthm}
\usepackage{booktabs}
\usepackage{array}
\usepackage{tabularx}
\usepackage{makecell}
\usepackage{multicol}
\usepackage{parskip}
\usepackage[section]{placeins}
\usepackage{float}
\usepackage{afterpage}
\usepackage{needspace}
\usepackage[table]{xcolor}
\usepackage{enumitem}
\usepackage{graphicx}
\usepackage{pdflscape}
\usepackage{tikz}
\usetikzlibrary{positioning,calc,arrows.meta,fit,backgrounds}
\usepackage[most]{tcolorbox}
\usepackage[numbers,sort&compress]{natbib}
\usepackage{caption}
\usepackage[colorlinks=true,linkcolor=blue!60!black,citecolor=blue!60!black,urlcolor=blue!60!black]{hyperref}

\usepackage{enumitem}

\definecolor{wwtdpurple}{RGB}{112,48,160}
\definecolor{SumBlue}{RGB}{33,82,140}
\definecolor{PosRed}{RGB}{176,58,46}
\definecolor{NoteTeal}{RGB}{0,122,118}
\definecolor{NoteViolet}{RGB}{104,52,160}
\setlist{nosep}

\newcommand{\norm}[1]{\lVert #1 \rVert}
\newif\ifshownotes \shownotestrue
\definecolor{NewGreen}{RGB}{0,110,40}
\newcommand{\vscl}{VSCL}
\newcommand{\veritas}{VERNA}

\newtcolorbox{summary}{breakable,enhanced,colback=SumBlue!6,colframe=SumBlue!6,
  borderline west={2.5pt}{0pt}{SumBlue},boxrule=0pt,frame hidden,
  left=8pt,right=8pt,top=5pt,bottom=5pt,arc=1pt,before skip=4pt,after skip=8pt,fontupper=\small}
\newtcolorbox{position}{enhanced,colback=PosRed!5,colframe=PosRed!5,
  borderline west={2.5pt}{0pt}{PosRed},boxrule=0pt,frame hidden,
  left=8pt,right=8pt,top=5pt,bottom=5pt,arc=1pt,before skip=6pt,after skip=6pt}

\newtcolorbox{example}{breakable,enhanced,colback=gray!6,colframe=gray!45,boxrule=0.4pt,arc=2pt,
  left=7pt,right=7pt,top=3pt,bottom=3pt,before skip=4pt,after skip=7pt,fontupper=\small}

\newtcolorbox{principlecard}[2]{breakable,enhanced,colback=white,colframe=SumBlue!30,boxrule=0.5pt,arc=2pt,
  borderline west={2.5pt}{0pt}{SumBlue},left=8pt,right=8pt,top=4pt,bottom=5pt,after skip=10pt,
  before={\par\needspace{7\baselineskip}\vspace{9pt}},
  title={\sffamily\bfseries #1\quad #2},colbacktitle=SumBlue!9,coltitle=SumBlue!70!black,
  toptitle=2.5pt,bottomtitle=2.5pt,titlerule=0pt,fonttitle=\small}

\DeclareUrlCommand\lpath{\urlstyle{tt}}   %

\title{%
Trustworthy Numerical Software in the Age of AI: \\What Would Turing Do?
}
\author{%
Kimon Fountoulakis\thanks{University of Waterloo} \and
Maksim Melnichenko\thanks{University of California, Berkeley (UCB)}\kern0.6em\thanks{International Computer Science Institute (ICSI)} \and
Muhammad Fetrat Qharabagh\footnotemark[1] \and
Georgios Alexandros Georgantzas\thanks{National Technical University of Athens (NTUA)} \and
George Spyrou\footnotemark[1] \and
Oleg Balabanov\footnotemark[2]\kern0.6em\footnotemark[3]\kern0.6em\thanks{Lawrence Berkeley National Laboratory (LBNL)} \and
Vasileios Georgiou \thanks{Purdue University}
 \and
 Riley Murray\thanks{Sandia National Laboratories} \and
Michael W.\ Mahoney\footnotemark[2]\kern0.6em\footnotemark[3]\kern0.6em\footnotemark[5]
}

\begin{document}
\maketitle

\begin{abstract}
AI coding agents can now produce numerical software faster than human experts can review it, raising serious concerns about usability and reliability. These and other AI tools are already transforming mathematics and software, and they will soon transform mathematical software. These concerns motivated us to assess AI-generated numerical linear algebra (NLA) kernels, which we did through an audit of a recent public GPU QR competition. All five leading kernels passed the competition’s accuracy tests, and the winning kernel achieved a 107-fold speedup over the competition’s baseline, which calls NVIDIA’s QR routine. Yet each failed at least once on new inputs, even those from the competition’s generator, under the same checker, while the baseline passed them all. In light of this, we want to develop reusable formalizations of numerical analysis methods, which remain sparsely represented in existing formal libraries. We describe VSCL, our Lean library of mathematical foundations for numerical computing. We also propose VERNA to connect these tasks by using formalized results to guide kernel assessment against explicit specifications. Drawing on Turing’s work on rounding-error analysis and program checking, we discuss design principles for developing trustworthy numerical software when direct human checking is impractical.
\end{abstract}

\definecolor{introrevisiongreen}{RGB}{0,100,0}
\providecommand{\introgreenrevision}[1]{%
  \begingroup
  \color{introrevisiongreen}%
  \ifdefined\hypersetup
    \hypersetup{allcolors=introrevisiongreen}%
  \fi
  #1%
  \endgroup
}

\section{Introduction}\label{sec:intro}

AI systems are increasingly used to construct mathematical arguments, formal
proofs, and numerical software~\cite{deepmind2025imo,tsoukalas2026alphaproofnexus,openai2026math,lanyon2026advdiff,lanyon2026maxwell}; and coding agents~\cite{jimenez2024swebench,yang2024sweagent} can now produce such software~\cite{demmel2026dbdsvdmr3,larsen2026lapackpr1377,larsen2026lapackpr1406,sakana2025,robustkbench2025} far faster than an expert can review it.
This creates challenges, since for over five decades the main assurance that numerical software can be trusted has been (human) expert authorship and (human) review~\cite[Section~2]{demmelgrading2026}.
When review cannot keep up, further evidence that numerical
software meets its stated requirements is needed. 
Here, we focus on numerical requirements,
such as accuracy within a specified tolerance~\cite[Section~10]{demmelgrading2026}, for numerical software, with the goal of assessing numerical linear algebra (NLA) kernels against stated accuracy requirements. 
We present an audit of NLA kernels from a recent public GPU competition; and, as we will describe, the results were sobering.
Motivated by these results, we aim to develop \vscl{}, a library of formalized numerical analysis consisting of reusable formalizations of numerical methods. 
We also propose the \veritas{} project, which aims to connect libraries of formalized numerical analysis and the assessment of numerical methods against stated accuracy requirements.  
It is well known that software can pass its evaluation and still fail in use~\cite{krentsel2026reality}, a problem that AI coding agents will only make worse. However, numerical software has unique properties, and numerical mathematics has a long history, with deep principles upon which we can draw.

\subsection{Assessing AI-generated NLA kernels: a cautionary example}

A timely illustration of reliability issues (that will become increasingly pressing) came in the summer of 2026, when the GPU MODE community~\cite{gpumode2026site} held four public competitions for fast NLA kernels~\cite{gpumode2026linalg,gpumode2026qr1,gpumode2026qr2,gpumode2026eigh,gpumode2026cholesky}.
The winning entry of the second QR competition ran 107 times faster than the competition's baseline~\cite{gpumode2026qr2,gpumode2026dataset,kernelbot2025}, which calls NVIDIA's QR routine through PyTorch~\cite{gpumode2026referencekernels}.
The organizers had invited submissions from ``you and your agents''~\cite{gpumode2026linalg}, and some participants reported building their entries with AI coding agents, one of them through more than 1{,}500 submissions in two weeks~\cite{sankalp2026autoresearch,lutz2026qr}.
Interested in the quality and performance of such AI-generated numerical software, we decided to audit the competition's five highest-ranked kernels (Section~\ref{sec:audit}).
All five had passed every accuracy test that the competition ran. 
We evaluated each kernel on 1,050 test cases comprising 945 distinct
inputs from the competition's own generator, using the competition's
own checker and keeping the kernels, accuracy requirements, and GPU
model unchanged~\cite{melnichenko2026auditrecord}.
The winning entry failed on 155 of these test cases, across nine
of the generator's ten matrix families; on some rank-deficient
inputs, it returned factors containing non-finite values.
Each of the other four entries also failed, two of them on about
a third of the test cases. The competition's baseline passed all
1,050 test cases.
Among other things, our audit demonstrated that passing a competition's original tests does \emph{not} ensure that the kernels can meet the same accuracy requirements on new inputs, even those generated from the same generator. 
This example illustrates the need for systematically assessing 
(in particular, AI-generated)
NLA kernels.

\subsection{Limited coverage of numerical analysis in formal libraries}

AI systems are used to develop and formalize mathematical results. Math Inc's Gauss
completed the formalization of the strong prime number theorem
\cite{mathinc2025strongpnt}; and Anthropic reported a Lean proof of
Fermat's Last Theorem that was produced largely autonomously in eleven days
\cite{anthropic2026flt}. 
On October 6, 2026, OpenAI released 722 machine-generated manuscripts, among
them a reported upper bound of $9/4$ on the exponent of matrix multiplication
over the complex numbers, with an accompanying Lean formalization
\cite{openai2026math,openai2026ninefourths}. On October 26, 2026, Axiom announced a complete
Lean formalization of the classification of finite simple groups,
with approximately 43.1 million lines of Lean. The 24-day effort
built on existing formalizations and used AI agents directed by
mathematicians and engineers \cite{axiom2026cfsg}. Such formalizations rely on libraries of definitions, theorems, and lemmas
checked by a proof assistant. Examples include Mathlib for Lean
\cite{mathlib2020cpp,mathlib4repo2026}, the Archive of Formal Proofs for Isabelle
\cite{afp2026}, Mathematical Components for Coq \cite{mahboubi2022mathcomp},
and the Mizar Mathematical Library \cite{bancerek2018mizar}.
Formal mathematical libraries could also provide foundations for software verification~\cite{rammal2026atlas,yang2026formalreasoning}.

\begin{figure}[htb!]
\centering
\includegraphics[width=\textwidth]{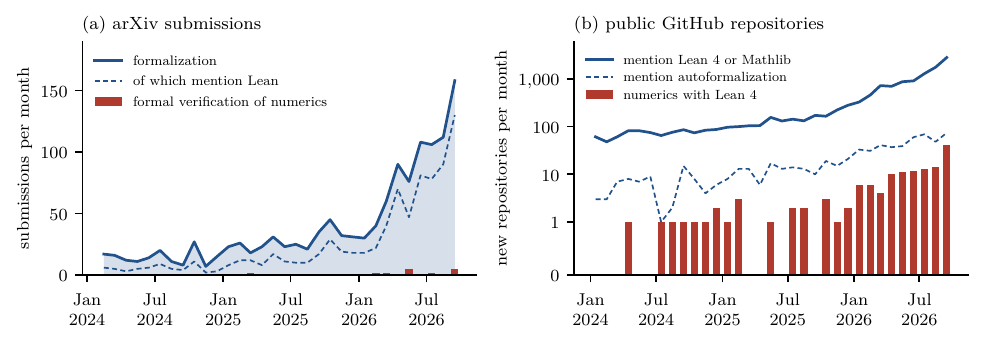}
\caption{Growth of formalization work, and how little of it concerns numerical analysis methods. (a) Monthly arXiv submissions, February 2024 to September 2026, whose title or abstract mentions autoformalization, or a form of the word formalization together with a proof assistant (Lean, Coq, Rocq, or Isabelle by name, or the phrases ``proof assistant'' and ``machine-checked'') (blue). The dashed line counts those that mention Lean or Mathlib, which account for most of the growth. Red bars count submissions on floating-point arithmetic or numerical analysis together with a proof assistant or formal verification. (b) Public GitHub repositories created each month, January 2024 to September 2026, whose name, description, or README mentions Lean~4 or Mathlib (solid), autoformalization (dashed), or Lean~4 together with floating point, rounding error, numerical analysis, or numerical linear algebra (red), on a scale that is logarithmic above one. The red bars total 140 repositories against 12{,}523 for the solid line, and 479 when Mathlib is counted as well. Part of the 2026 rise reflects repositories generated by coding agents. Both panels were retrieved through the public arXiv and GitHub search interfaces on November 11, 2026.}
\label{fig:activity}
\end{figure}

The growth of formalized mathematics is also visible in research papers and in software.
ArXiv submissions matching our queries for formalization
rose from about fifteen per month in early 2024 to about 160 per month in
September 2026 (Figure~\ref{fig:activity}a). GitHub repositories
mentioning Lean~4 or Mathlib went from about sixty a month to almost
three thousand in that time (Figure~\ref{fig:activity}b).

\begin{figure}[htb!]
\centering
\includegraphics[width=\textwidth]{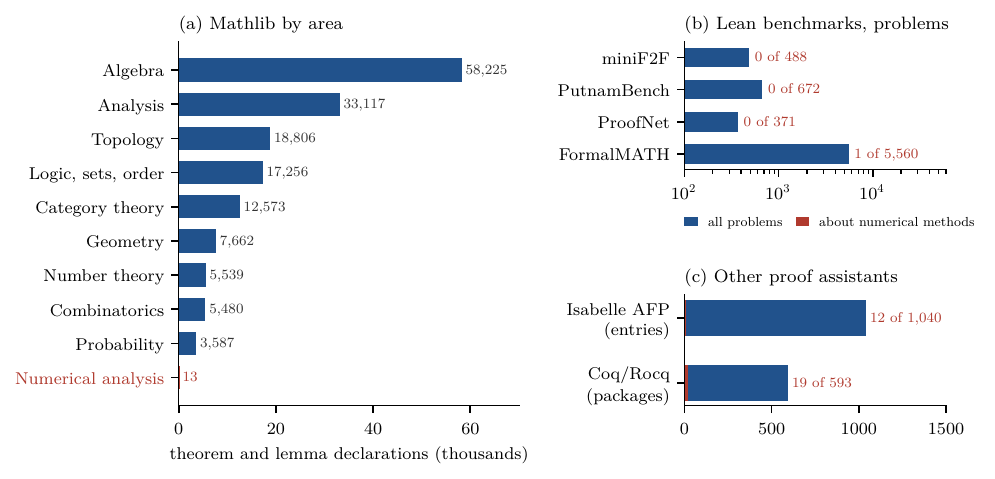}
\caption{How little formal mathematics concerns numerical methods. (a) Lines that begin a theorem or lemma in Mathlib, by area, at commit \texttt{021ce68} of October 7, 2026. A keyword search finds 13 for numerical analysis in two files: an experimental floating-point model and the error bound of the trapezoidal rule. (b) Problems in four Lean benchmarks, miniF2F \cite{minif2f2022}, PutnamBench \cite{putnambench2024}, ProofNet \cite{proofnet2023}, and FormalMATH \cite{formalmath2025}, and how many concern a numerical method. The one FormalMATH problem runs three steps of gradient descent in exact arithmetic. (c) Entries in Isabelle's Archive of Formal Proofs and released Coq and Rocq packages, and how many concern numerical methods, mostly floating-point and interval arithmetic. Counts in (a) to (c) come from keyword searches whose hits we read one by one, and the units differ by panel.}
\label{fig:subfields}
\end{figure}

Despite the growing popularity of formalized mathematics, numerical analysis remains sparsely represented in the formal libraries and benchmarks we surveyed.
The gap is visible in Mathlib \cite{mathlib4repo2026}, the Lean library on
which most of the systems above build.
A keyword-based survey of its roughly 182{,}000 theorem and lemma
declarations identifies 13 that concern finite precision, conditioning, stability, or convergence in numerical methods\footnote{Mathlib holds 104 declarations on exact-arithmetic matrix algorithms. These declarations are excluded from our numerical-analysis count because none involves finite precision, conditioning, stability, or convergence.}.
All 13 are in just two files, one on an experimental
floating-point model and the other on the error bound for the trapezoidal
rule. Numerical methods are also rare in Lean benchmarks. Across four
widely used benchmarks, only one of the 7{,}091 problems concerns a
numerical method (Figure~\ref{fig:subfields}). Beyond Lean, 12 of the
1{,}040 entries in Isabelle's Archive of Formal Proofs
\cite{afp2026submission} and 19 of the 593 released Coq and Rocq packages
concern numerical methods, mostly floating-point and interval arithmetic,
and we found no entry in the Archive on the rounding-error analysis of
linear algebra (Figure~\ref{fig:subfields}c).
We refer the reader to Appendix~\ref{sec:prior} for a more detailed literature review on this subject.

\subsection{Our contributions and research direction}\label{sec:intro-contributions}

This paper makes several contributions: we describe an audit of NLA kernels from a recent public GPU competition; we present the first version of a library of formalized numerical analysis; we propose a project to connect libraries of formalized numerical analysis with the assessment of numerical methods against stated accuracy requirements; and we describe principles that should guide efforts such as these, which involve formalizing numerical mathematics.
\begin{itemize}
\item
First, we present the QR kernel audit described above; see Section~\ref{sec:audit}. 
Its failures are recorded and can be replayed~\cite{melnichenko2026auditrecord}; and they highlight the need for efforts such as what we present. 
\item
Second, we describe our Verified Scientific Computing Library (\vscl{}) project, which provides reusable formalized numerical analysis; see Section~\ref{sec:vscl}.
\vscl{} contains selected definitions, algorithms, and theorems in numerical analysis and its foundations, with proofs checked by Lean~\cite{lcm2026repo}.
Its sources include Higham's \emph{Accuracy and Stability of Numerical Algorithms}~\cite{higham2002asna}, LeVeque's \emph{Finite Volume Methods for Hyperbolic Problems}~\cite{leveque2002fvm}, and Vershynin's \emph{High-Dimensional Probability}~\cite{vershynin2018hdp}.
It also includes selected material from Eymard, Gallou{\"e}t, and Herbin's \emph{Finite Volume Methods}~\cite{eymard2019fvm} and Nocedal and Wright's \emph{Numerical Optimization}~\cite{nocedal2006numerical}.
The coverage includes rounding-error and perturbation analysis, matrix algorithms, finite-volume methods, probability, and optimization.
\item
Third, we propose to connect these efforts through the Verification Engine for Reliable Numerical Algorithms (\veritas{}); see Section~\ref{sec:veritas}.
\veritas{} is designed to assess NLA kernels against explicit specifications, such as a kernel's accuracy requirements, using formalized results from \vscl{} to guide the choice
of checks and inputs. 
This paper states the goal of \veritas{} and leaves its
design description to a separate publication.
\item 
Fourth, we discuss what lies ahead, including fundamental principles that should guide the development of formalized numerical mathematics; see Section~\ref{sec:conclusion}.
We arrived at these by asking what Turing, one of the founders of computer science~\cite{turing1936computable} who also made fundamental conceptual contributions to early numerical analysis~\cite{turing1948rounding} and program checking~\cite{turing1949checking}, would do if he were faced with the AI-generated challenges we are facing today.
\end{itemize}

Our website~\cite{tnpg2026website} provides links to the audit results, a catalogue of the results formalized in \vscl{}, and the statements and proofs for the case study in Section~\ref{sec:randnla}. These materials are publicly available, with links to their sources, so readers can check our claims directly.

The rest of this paper is organized as follows. Section~\ref{sec:audit} presents our audit of the GPU MODE competition.
Section~\ref{sec:vscl} presents \vscl{}, its mathematical content,
the reuse study, and the formalization workflow. It also applies the library to a survey of randomized numerical linear algebra (Section~\ref{sec:randnla}). Section~\ref{sec:veritas} introduces \veritas{} and its three steps.
Section~\ref{sec:conclusion} discusses the remaining work on underflow and
exceptional values and the cost of formalizing at scale; and it closes with the
Turing-inspired goals for the broader research program.
Appendix~\ref{sec:prior} surveys prior work on formalized numerical analysis and kernel verification, and Appendices~\ref{sec:audit-details}, \ref{app:benchmark_details} and~\ref{sec:open-weights} give the details of the audit, the reuse analysis, and our experiments with open-weight models.

\section{An audit of the GPU MODE competition}
\label{sec:audit}

GPU MODE is a community of GPU developers whose competitions reward fast implementations of specified operations~\cite{gpumode2026retrospective,gpumode2026site,gpumode2026linalg}.
In the summer of 2026, it held four linear algebra competitions: two for QR, one for the symmetric eigendecomposition, and one for Cholesky~\cite{gpumode2026qr1,gpumode2026qr2,gpumode2026eigh,gpumode2026cholesky}.\footnote{The organizers announced one competition with three problems~\cite{gpumode2026linalg}, and QR ran on two leaderboards. We refer to each of the four leaderboards as a competition.}
We performed an audit of all four competitions.

This section reports our audit of the second GPU MODE QR competition, which asked participants to submit kernels that factor batches of single-precision square matrices.
(Appendix~\ref{sec:audit-details} reports our audits of the other three competitions.)
Each entry was first checked on 22 test cases and then timed and checked on 12 additional test cases, which we call the \textit{scored cases}.
Entries that passed every check were ranked by speed on NVIDIA B200 GPUs, with compute funded by Core Automation~\cite{gpumode2026linalg}.
The baseline was the competition's reference implementation, which calls NVIDIA's QR routine through PyTorch's \texttt{torch.geqrf}~\cite{paszke2019pytorch,gpumode2026referencekernels,nvidiacusolver2026}.

We first became interested in the competition because of the winning entry's speed: it ran \textit{107 times faster than the baseline} (Appendix~\ref{sec:audit-env})~\cite{gpumode2026qr2,gpumode2026dataset,kernelbot2025}.
Our interest grew when participants reported building their entries with AI coding agents.
One participant describes more than 1{,}500 submissions made with a coding agent in two weeks, and a twelfth place among 183 participants~\cite{sankalp2026autoresearch}.
A participant whose entry we audited describes nearly 2{,}000 experiments run with coding agents in fifteen days~\cite{lutz2026qr}.
Both participants stated that QR was new to them (which also intrigued us).
The speed and these accounts led us to audit the five highest-ranked entries on a wider range of inputs that the competition allowed.

\subsection{QR accuracy requirements}

A QR factorization writes a matrix $A$ as the product of an orthogonal matrix $Q$ and an upper triangular matrix $R$.
A kernel computes these factors in floating-point arithmetic, so they satisfy this definition approximately.
The competition's checker requires every returned value to be finite and both residuals in Table~\ref{tab:qr-contract} to be within their tolerances.
The first residual measures how well the computed factors reproduce $A$, and the second measures how far the computed $Q$ is from orthogonal.

\begin{table}[!htbp]
\centering
\small
\caption{The competition's accuracy requirements for the QR factorization of a square matrix $A$ of order $n$.
Here $\norm{\cdot}_1$ is the matrix 1-norm and $\varepsilon=2^{-23}$ is single-precision machine~epsilon.}
\label{tab:qr-contract}
\renewcommand{\arraystretch}{1.15}
\begin{tabular}{@{}lcc@{}}
\toprule
Property & Residual & Tolerance \\
\midrule
Factorization & $\norm{R-Q^T A}_1$ & $20n\varepsilon\norm{A}_1$ \\
Orthogonality & $\norm{Q^TQ-I}_1$ & $100n\varepsilon$ \\
\bottomrule
\end{tabular}

\end{table}

The organizers scaled both tolerances with the order $n$, as is standard~\cite{gpumode2026linalg}.
The checker applies both requirements even when the input is rank deficient or poorly conditioned.
A stable Householder QR algorithm handles such inputs while keeping both residuals small.

\subsection{Exposing the failures}

Our audit consists of 1,050 test cases comprising 945 distinct
inputs. These test cases were generated using the organizers'
generator across ten matrix families, five conditioning settings,
and three seeds for each shape~\cite{gpumode2026referencekernels}.
All five highest-ranked entries failed on some of the inputs of this audit.
Results are summarized in Table~\ref{tab:qr-second-round}, including two reference implementations.
The bfloat16-rounded baseline provides a sense of what might happen with a very particular (and very dubious) attempt at optimizing kernel performance.

\begin{table}[!htbp]
\centering
\small
\caption{QR entries and controls (Appendix~\ref{sec:audit-qr-inputs}) under the competition's checker.
Scored failures: rejected outputs on the twelve scored cases.
Audit failures: rejected outputs in the 1,050 audit evaluations,
covering 945 distinct inputs.
Families with failures: matrix families, out of ten, with at least one audit failure.
Data: audit record~\cite{melnichenko2026auditrecord}.
}
\label{tab:qr-second-round}
\renewcommand{\arraystretch}{1.15}
\begin{tabular}{@{}clrrr@{}}
\toprule
Leaderboard rank & Entry & Scored failures & Audit failures & Families with failures \\
\midrule
1 & 10billiontokens & 0/12 & 155/1{,}050 & 9/10 \\
2 & gau.nernst & 0/12 & 46/1{,}050 & 2/10 \\
3 & dhu.randhar & 0/12 & 327/1{,}050 & 9/10 \\
4 & zhongmingee & 0/12 & 3/1{,}050 & 1/10 \\
5 & michaelmelons & 0/12 & 357/1{,}050 & 10/10 \\
\midrule
\multicolumn{2}{@{}l}{Baseline} & 0/12 & 0/1{,}050 & 0/10 \\
\multicolumn{2}{@{}l}{Baseline, bfloat16-rounded input} & 8/12 & 619/1{,}050 & 10/10 \\
\bottomrule
\end{tabular}

\end{table}

\begin{figure}[!htbp]
\centering
\includegraphics[width=\linewidth]{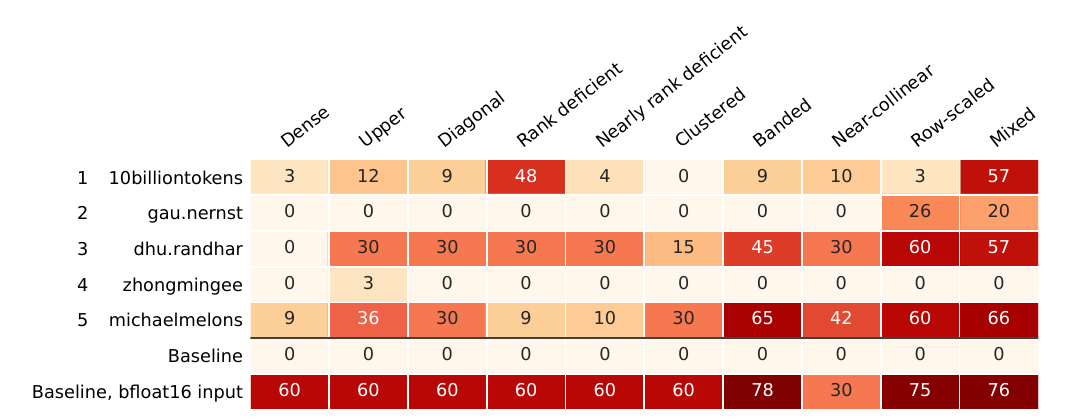}
\caption{Audit failures by matrix family.
Rows: the five audited entries by leaderboard rank, then the two controls.
Each cell counts rejected outputs among the 105 audit inputs of a family.
Data: audit record~\cite{melnichenko2026auditrecord}.}
\label{fig:qr-families}
\end{figure}

Different entries failed on different matrix families (Figure~\ref{fig:qr-families}).
For the winning entry, every rank-deficient input of order 32, 176, or 352 produced non-finite values in the returned factors.
The runner-up's 46 failures were confined to row-scaled matrices and batches mixing several conditioning profiles.
Four of the five entries also test their input for the generator's constructions at scored orders and run dedicated code for them~\cite{gpumode2026dataset}.
One inspects four matrices of a batch when the order is 1024.
When they match the generator's nearly rank-deficient construction, it factors only 768 of the 1{,}024 columns.

\subsection{Audit conclusions}

Given that the entries passed every check that the competition ran, our main conclusions concern what ``passing a fixed set of checks'' can really establish for software verification.

\begin{enumerate}
\item \textbf{Trust needs a stated accuracy requirement, met on every code path.}
The competition's checker states such a requirement, the two tolerances of
Table~\ref{tab:qr-contract}, and a kernel must meet it on every code path.
The fourth-ranked competition entry used its own fast QR methods for some inputs and a standard QR routine as a fallback for others. In our audit, every output identified as a fallback passed the checker, while its own methods failed on three inputs.

\item \textbf{Passing a fixed test set makes no statement about the general reliability of the kernel.}
All five entries passed the twelve scored cases and failed on other inputs from the same generator.
Tests exposed these failures.
Tests cannot show that no more remain.

\item \textbf{Fixed tests can be gamed.}
Four of the five entries recognize the generator's constructions and handle them with dedicated code.
Demmel and co-authors consider aggressive gaming more likely from an AI than from a human developer~\cite[Section~4]{demmelgrading2026}.
They design their tests to make it expensive and state that they have no proof that this suffices~\cite[Sections~2 and~4]{demmelgrading2026}.

\item \textbf{Authorship alone does not establish reliability.}
The organizers invited entries from ``you and your agents''~\cite{gpumode2026linalg}.
Two participants who described themselves as new to QR reported building high-ranking entries with AI coding agents in about two weeks~\cite{sankalp2026autoresearch,lutz2026qr}.
One of the two entries ranked fifth, and so was among the five highest-ranked entries that we audited.
All five failed on new inputs, the fifth-ranked one in all ten matrix families (Table~\ref{tab:qr-second-round}).
Numerical code can no longer be assumed to come from trustworthy teams of professional programmers~\cite[Section~2]{demmelgrading2026}.

\end{enumerate}

Our audit's conclusions call for an assessment of NLA kernels that does not depend (only) on a fixed test set or on trust in a kernel's author.
Regarding software verification more generally, Krentsel and co-authors study why software can pass its evaluator and still fail in use, a problem that AI coding agents make worse~\cite{krentsel2026reality}.
They note that this problem becomes more tractable when ``\textit{intent is fixed by an authoritative artifact}'', when ``\textit{the relevant world is bounded by a stable formal interface}'', or when ``\textit{an existing implementation provides a strong behavioral reference}''~\cite[Appendix~A]{krentsel2026reality}.
We argue that NLA meets these three conditions.\footnote{Krentsel and co-authors name a fourth condition, that ``\textit{a faithful model of the world can be learned from available data}''. NLA has no need of it, because the world of a kernel is the arithmetic of its hardware, which is specified and does not have to be learned.}
Error bounds fix the largest error a kernel may produce, floating-point standards specify the arithmetic in which it runs and its exceptional values, and standard libraries such as BLAS and LAPACK provide the reference.
Our \vscl{} project (Section~\ref{sec:vscl}) formalizes error bounds in Lean and aims to connect them to the arithmetic that floating-point standards specify.
Unlike a test set, an error bound covers every input that satisfies its assumptions.

\section{\vscl{}}
\label{sec:vscl}

\vscl{} is a library that contains Lean definitions, algorithms, and proved results in numerical analysis and its mathematical foundations. 
It covers rounding-error and perturbation analysis, matrix algorithms, finite-volume methods, probability, and optimization. 

\subsection{The library}
\label{sec:library}

\vscl{} contains definitions, algorithms, and proved results that
other formalizations can use. It is publicly available as
\texttt{lean-computational-mathematics}~\cite{lcm2026repo}.
The coverage and size reported below refer to the published November 9,
2026 code snapshot. The development is an ongoing effort.

\paragraph{The mathematics.} 
The library contains selected material
from all twenty-eight chapters of Higham's \emph{Accuracy and Stability of
Numerical Algorithms}~\cite{higham2002asna}. This covers floating-point
arithmetic, the accumulation of rounding error, matrix norms and perturbation
theory, and the analysis of algorithms for linear systems and least-squares
problems. It also contains selected material from Chapters~1--4 of both LeVeque's
\emph{Finite Volume Methods for Hyperbolic Problems}~\cite{leveque2002fvm}
and Vershynin's \emph{High-Dimensional Probability}~\cite{vershynin2018hdp}.
The LeVeque material covers conservation laws, advection and wave
propagation, and the analysis of finite-volume methods, including the
Chapter~4 results on convergence, accuracy, and stability described below.
The Vershynin material covers random variables, sums of independent random
variables, concentration inequalities, and the Chapter~4 random-matrix
results on singular values, covariance estimation, and spectral clustering.
Recent additions include finite-volume foundations from Eymard, Gallou\"et
and Herbin~\cite{eymard2019fvm},
and optimization foundations from Nocedal and Wright~\cite{nocedal2006numerical}.
These developments extend the library beyond rounding-error bounds to
higher level numerical algorithms and the probability theory used to
analyze randomized algorithms.
The coverage of each
source is given below.

\textit{Accuracy and Stability of Numerical Algorithms}~\cite{higham2002asna}. The library covers claims
across all twenty-eight chapters
(Table~\ref{fig:coverage-higham}). 
We distinguish claims
proved as stated, printed claims with formal counterexamples and
separately proved corrections, and deferred items whose source is
imprecise or refers to external work~\cite{lcm2026repo}.

\begingroup
\let\vsclOriginalCaption\caption
\renewcommand{\caption}{%
  \par\vspace{3pt}%
  \vsclOriginalCaption}
\begin{table}[!htb]
\centering
\begingroup
\normalfont\small
\setlength{\parindent}{0pt}
\setlength{\parskip}{0pt}
\newcommand{\covrow}[2]{%
  \noindent\begin{minipage}[t]{0.095\linewidth}\raggedright
  #1\strut\par
  \end{minipage}\hfill
  \begin{minipage}[t]{0.875\linewidth}\raggedright
  #2\strut\par
  \end{minipage}\par\vspace{4pt}}
\begin{minipage}{\linewidth}
\textbf{Accuracy and Stability of Numerical Algorithms~\cite{higham2002asna}}\par
\vspace{4pt}\hrule height 0.6pt\vspace{5pt}
\covrow{Ch.\ 1--5}{\emph{Finite precision and basic computations.} Floating-point principles and arithmetic; dot products and accumulated rounding error; summation; polynomial evaluation and Horner's method.}
\covrow{Ch.\ 6--8}{\emph{Norms and triangular systems.} Vector and matrix norms; perturbation theory for linear systems; triangular solves.}
\covrow{Ch.\ 9--13}{\emph{Factorizations.} LU, Cholesky, and block LU; symmetric indefinite and skew-symmetric systems; iterative refinement.}
\covrow{Ch.\ 14--18}{\emph{Matrix computations.} Inversion; condition-number estimation; the Sylvester equation; stationary iterations; matrix powers.}
\covrow{Ch.\ 19--22}{\emph{QR and least squares.} QR factorization; least-squares problems; underdetermined and Vandermonde systems.}
\covrow{Ch.\ 23--28}{\emph{Further methods.} Fast matrix multiplication; fast Fourier transforms and applications; nonlinear systems and Newton's method; automatic error analysis; floating-point software; test matrices.}
\hrule height 0.4pt
\end{minipage}
\endgroup
\caption{Material across all 28 chapters of~\cite{higham2002asna}. Claims proved as stated, corrected discrepancies, and deferred statements are recorded separately~\cite{higham2002asna,lcm2026repo}.}
\label{fig:coverage-higham}
\end{table}

\endgroup

\textit{Finite Volume Methods for Hyperbolic Problems}~\cite{leveque2002fvm}.
The library contains material from Chapters~1--4
(Table~\ref{fig:coverage-leveque}). This includes integral and
differential formulations of conservation laws, scalar advection and
linear acoustics, and the cell averages and flux updates used in
finite-volume methods. It also covers characteristic curves, wave
superposition, domains of dependence and influence, and linear Riemann
solutions, together with material on boundary and periodic problems.
Chapter~4 adds convergence and Courant--Friedrichs--Lewy (CFL) necessity
results for specified data classes, norms, and refinements, together with
accuracy and stability results, explicit and Crank--Nicolson diffusion,
Godunov methods, and linear wave propagation.

\begin{table}[!htb]
\centering
\begin{minipage}{\linewidth}
\small
\noindent\textbf{Finite Volume Methods for Hyperbolic Problems}~\cite{leveque2002fvm}\par
\vspace{3pt}\hrule height 0.6pt\vspace{5pt}
\renewcommand{\arraystretch}{1}
\begin{tabular}{@{}p{0.15\linewidth}@{}p{0.85\linewidth}@{}}
Chapter 1 & \raggedright\emph{Conservation laws and finite-volume ideas.}
Fluxes, scalar advection, linear acoustics, hyperbolicity and wave speeds;
cell averages and flux updates; interface Riemann data and dimensional
splitting.\tabularnewline[0.5em]
Chapter 2 & \raggedright\emph{Physical equations.}
Integral and differential conservation; advection, diffusion and reaction;
fluids, acoustics, elasticity and waves; Eulerian and Lagrangian coordinates;
linearization and hyperbolicity.\tabularnewline[0.5em]
Chapter 3 & \raggedright\emph{Characteristics and Riemann problems.}
Characteristic coordinates and curves; wave superposition and initial data;
domains of dependence and influence; linear Riemann solutions, acoustic
waves and strengths, boundary and periodic problems.\tabularnewline[0.5em]
Chapter 4 & \raggedright\emph{Finite-volume methods for linear equations.}
Numerical domains of dependence; convergence and CFL necessity under
stated assumptions; accuracy and stability; explicit and Crank--Nicolson
diffusion; reconstruction, evolution and averaging; Godunov methods and
linear wave propagation.\tabularnewline
\end{tabular}
\par\vspace{4pt}\hrule height 0.4pt\vspace{3pt}
\end{minipage}
\caption{Material from Chapters~1--4 of~\cite{leveque2002fvm}
in the November 9, 2026 snapshot.}
\label{fig:coverage-leveque}
\end{table}

\textit{High-Dimensional Probability}~\cite{vershynin2018hdp}.
The library contains selected material from Chapters~1--4
(Table~\ref{fig:coverage-vershynin}).
This includes distributions, expectation, moments and tail bounds,
together with laws of large numbers and Poisson limits. It also covers
sums of independent random variables, sub-Gaussian and sub-exponential
variables, Orlicz norms, and Hoeffding, Chernoff, Bernstein, and
bounded-difference inequalities. The Chapter~3 material
includes random vectors, isotropy and covariance, as well as proofs
of Euclidean norm concentration, Grothendieck's inequality, and the
Max-Cut approximation result.
The Chapter~4 development adds singular-value bounds, covering and packing,
coding, random-matrix concentration, covariance estimation, and spectral
clustering.

\begin{table}[!htb]
\centering
\begin{minipage}{\linewidth}
\small
\noindent\textbf{Vershynin} \emph{High-Dimensional Probability}\par
\vspace{3pt}\hrule height 0.6pt\vspace{5pt}
\renewcommand{\arraystretch}{1}
\begin{tabular}{@{}p{0.15\linewidth}@{}p{0.85\linewidth}@{}}
Chapter 1 & \raggedright\emph{Random variables.}
Distributions, expectation, variance, moments, tails and moment-generating
functions; probability inequalities; laws of large numbers and Poisson limits.
\tabularnewline[0.5em]
Chapter 2 & \raggedright\emph{Independent sums.}
Coin sums, Hoeffding and Chernoff bounds, and random-graph degrees;
sub-Gaussian and sub-exponential variables and sums; Orlicz norms,
Bernstein and bounded-difference inequalities.\tabularnewline[0.5em]
Chapter 3 & \raggedright\emph{Random vectors.}
Norm concentration, isotropy, covariance and principal components;
high-dimensional distributions and sub-Gaussian vectors; Grothendieck's
inequality, semidefinite programs, Max-Cut and randomized rounding;
kernels and tensor products.\tabularnewline[0.5em]
Chapter 4 & \raggedright\emph{Random matrices.}
Singular values and operator norms; covering, packing and nets;
error-correcting codes; random-matrix concentration; covariance
estimation and spectral clustering.\tabularnewline
\end{tabular}
\par\vspace{4pt}\hrule height 0.4pt\vspace{3pt}
\end{minipage}
\caption{Material from Chapters~1--4 of~\cite{vershynin2018hdp}
in the November 9, 2026 snapshot.}
\label{fig:coverage-vershynin}
\end{table}

\textit{Finite Volume Methods}~\cite{eymard2019fvm}. The Chapter~1 development of
\emph{Finite Volume Methods} contains three
results concerning local conservativity, cell unknowns,
and conservative exchange.

\textit{Numerical Optimization}~\cite{nocedal2006numerical}. The Chapter~1 development of
\emph{Numerical Optimization} covers constrained
minimization, transportation models, finite candidate domains, and
convexity.

\paragraph{Size.} \vscl{} contains 5{,}003 production
Lean modules and 1{,}700{,}804 nonblank lines of production source:
1{,}336 reusable modules, 3{,}197 source modules, 460 aggregates,
five internal modules, and five attributed upstream modules.
A count across all 5{,}190 Lean files, including tests and tooling, finds
48{,}122 lines that begin a theorem or lemma declaration, with comments
excluded.

\paragraph{Reusing results.}
We group shared definitions and theorems by mathematical subject,
following Mathlib's conventions~\cite{mathlib4repo2026}. Separate
files for each book record which textbook claims these results
establish and explain any differences in notation or formulation.
When proving a new result, we can import the definitions and
theorems from the relevant subject files. The files for
each book preserve the connection to the source text, but need not
be imported to use the shared results~\cite{lcm2026repo}.

\textit{The library reuses its own results internally}. A dependency audit
of Higham's \emph{Accuracy and Stability of
Numerical Algorithms}~\cite{higham2002asna} examined references in Lean's elaborated statements
and proofs. Of 47{,}890 public declarations
identified as source-written, 35{,}889 (74.9\%) are referenced by another
project declaration. In addition, 46{,}811 (97.7\%) reach a dependency or
consumer through at least two successive references
(Table~\ref{tab:library-reuse})~\cite{lcm2026repo}.
For example, the matrix-multiplication error proof uses the
matrix-vector error result, which uses the dot-product bound,
which in turn uses the summation analysis and the addition-rounding
model. A longer QR-solve path
contains fourteen direct proof dependencies across seven modules,
from the solver's backward-error bound to the square-root rounding
model~\cite{lcm2026repo}. These references demonstrate formal composition within the library.

\providecommand{\BenchTimeGain}{31.5\%}
\providecommand{\BenchTimeWins}{38}
\providecommand{\BenchCatOneTimeGain}{-2.6\%}
\providecommand{\BenchCatOneTimeLoss}{2.6\%}
\providecommand{\BenchCatTwoTimeGain}{56.8\%}
\providecommand{\BenchCatTwoTimeLoss}{56.8\%}
\providecommand{\BenchCatThreeTimeGain}{61.8\%}
\providecommand{\BenchCatThreeTimeLoss}{61.8\%}
\providecommand{\BenchZeroTimeGain}{4.3\%}
\providecommand{\BenchLastTimeGain}{65.7\%}
\providecommand{\BenchTokenGain}{22.1\%}
\providecommand{\BenchTokenWins}{34}
\providecommand{\BenchCatOneTokenGain}{-10.5\%}
\providecommand{\BenchCatOneTokenLoss}{10.5\%}
\providecommand{\BenchCatTwoTokenGain}{38.8\%}
\providecommand{\BenchCatTwoTokenLoss}{38.8\%}
\providecommand{\BenchCatThreeTokenGain}{56.7\%}
\providecommand{\BenchCatThreeTokenLoss}{56.7\%}
\providecommand{\BenchZeroTokenGain}{23.4\%}
\providecommand{\BenchLastTokenGain}{51.5\%}
\providecommand{\BenchCodeGain}{29.3\%}
\providecommand{\BenchCodeWins}{42}
\providecommand{\BenchCatOneCodeGain}{-1.8\%}
\providecommand{\BenchCatOneCodeLoss}{1.8\%}
\providecommand{\BenchCatTwoCodeGain}{41.7\%}
\providecommand{\BenchCatTwoCodeLoss}{41.7\%}
\providecommand{\BenchCatThreeCodeGain}{59.1\%}
\providecommand{\BenchCatThreeCodeLoss}{59.1\%}
\providecommand{\BenchZeroCodeGain}{3.9\%}
\providecommand{\BenchLastCodeGain}{54.2\%}
\providecommand{\BenchUsingTasks}{33}
\providecommand{\BenchZeroTasks}{27}

\textit{Library reuse in new tasks can outweigh access overhead.}
The library can also help an agent tasked with proving a result complete the proof more efficiently when parts of the proof have already been formalized in the library. Although this benefit might seem obvious at first glance, providing access to a library incurs overhead, particularly for an agent that has not been specifically trained on its contents. This overhead includes searching the library for potentially useful results, assessing whether the results found are relevant to the task, and determining whether they can be used in the proof. These costs could outweigh the time and tokens saved by reusing existing results instead of proving them from scratch. We show, however, that when the library contains results supporting parts of a proof, the savings in time and tokens more than offset this overhead. Details of this experiment are given in Appendix~\ref{app:benchmark_details}. We briefly discuss the experiment below.

To test this, we extracted 60 statements from 54 research papers covering topics related to the chapters of Higham's \emph{Accuracy and Stability of Numerical Algorithms}~\cite{higham2002asna}. 
We evaluated the same agent on these statements under two conditions: with and without read-only access to \vscl{}, repeating each task three times per condition. All 360 attempts produced accepted proofs. However, our measurements show that reusing results from \vscl{} rather than proving them from scratch significantly reduces execution time, token consumption, and the number of lines of submitted Lean code, more than offsetting the overhead of library access. 
This finding applies to tasks whose proofs can use one or more library results that meet three criteria: (i) reconstructing the result without access to \vscl ~would require significant mathematical effort, as judged by the LLM that extracted the task (OpenAI Codex with GPT-6 Sol~\cite{openai2026codex}); (ii) the result goes beyond basic algebraic manipulations; and (iii) Mathlib does not already provide a readily usable equivalent. 

We grouped the tasks into three categories according to the number of \vscl{} results judged to satisfy the three conditions applicable to their proofs\footnote{These criteria ensure that tasks are not classified based on trivial results, such as elementary algebraic identities. For example, \vscl{} contains an established identity stating that subtraction distributes over finite expectations. Without these criteria, the extracting agent might classify a task as benefiting from library reuse simply because such an elementary result could be used as a minor step in its proof. This would distort the task classification and undermine the practical significance of the subsequent measurements and conclusions.}. The first category (C1) contains 30 tasks with no such supporting result in the library, the second (C2) contains 20 tasks with exactly one, and the third (C3) contains 10 tasks with at least two such results in the library.

For C2, access to \vscl{} reduced the category means for execution time, total token consumption, and lines of submitted Lean code by \BenchCatTwoTimeGain, \BenchCatTwoTokenGain, and \BenchCatTwoCodeGain, respectively. For C3, the corresponding reductions were \BenchCatThreeTimeGain, \BenchCatThreeTokenGain, and \BenchCatThreeCodeGain. For C1, access increased these three costs by \BenchCatOneTimeLoss, \BenchCatOneTokenLoss, and \BenchCatOneCodeLoss, respectively. These increases are consistent with the overhead of accessing and searching the library when no supporting result satisfying the three criteria is available to offset those costs. Appendix~\ref{app:benchmark_details} reports the analysis in full.

\begin{table}[h!]
\centering
\small
\begingroup
\renewcommand{\arraystretch}{1.12}
\begin{tabularx}{\linewidth}{@{}Xcc@{}}
\toprule
Use within the library & Declarations & Percentage \\
\midrule
Used by another declaration & 35{,}889 & 74.9\% \\
Linked through two or more steps & 46{,}811 & 97.7\% \\
No links to other declarations & 143 & 0.3\% \\
\bottomrule
\end{tabularx}
\endgroup
\caption{Internal reuse of formalized results in VSCL for Higham's \emph{Accuracy and Stability of
Numerical Algorithms}~\cite{higham2002asna}. The audit counts 47{,}890
public declarations identified in the source code, such as definitions and
theorems. A link means that one declaration uses another in its statement or
proof. Two steps involve using a declaration that itself uses another, or being
used through such a chain. The rows overlap, so the percentages should
not be added~\cite{lcm2026repo}.}
\label{tab:library-reuse}
\end{table}

\subsection{The formalization workflow}
\label{sec:skill}

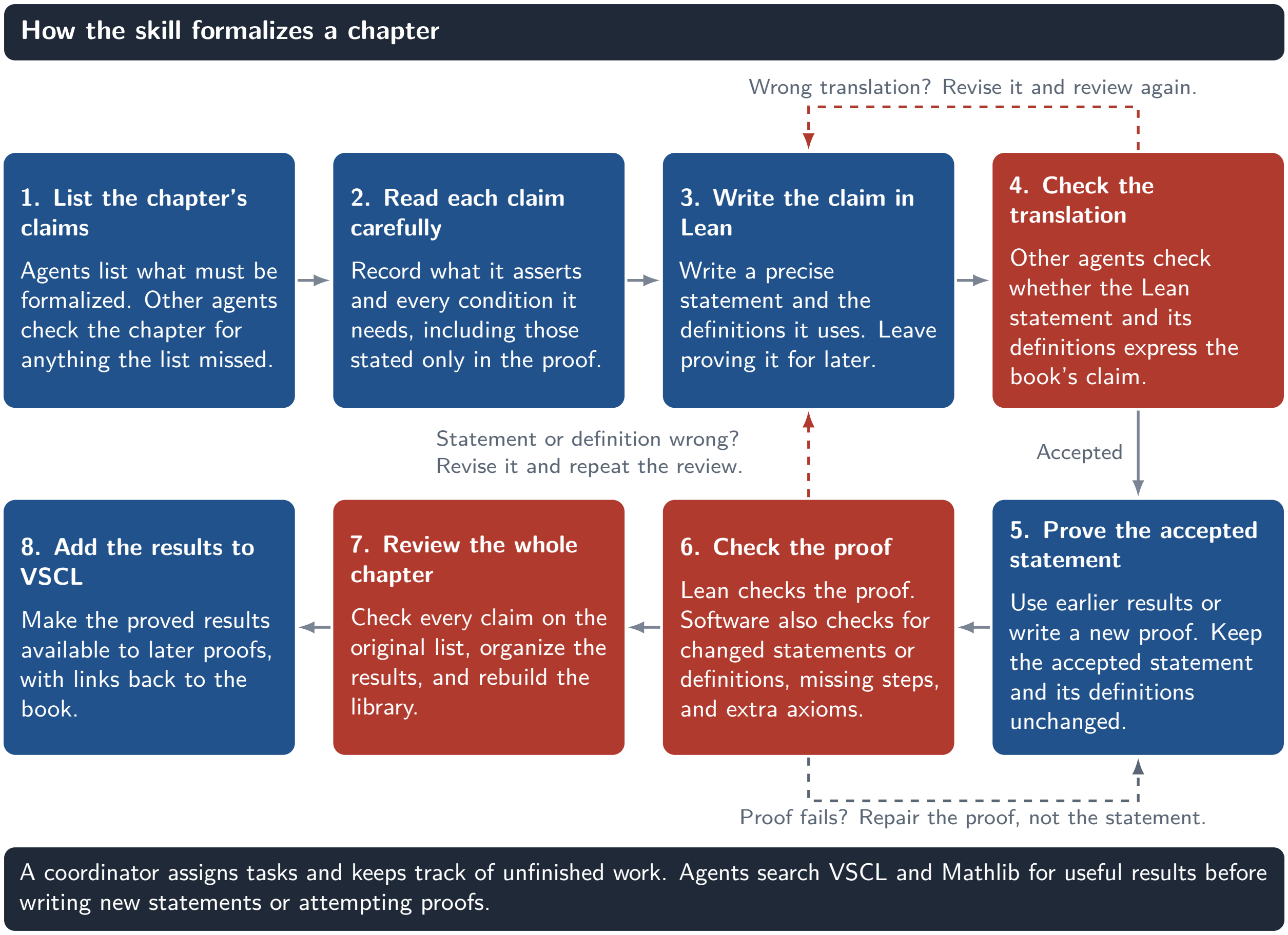
\begin{figure}[!htb]
\centering
\definecolor{VSCLInk}{HTML}{1F2937}
\definecolor{VSCLArrow}{HTML}{7A8492}
\definecolor{VSCLMuted}{HTML}{5B6675}
\resizebox{\linewidth}{!}{%
\begin{tikzpicture}[
  font=\sffamily\fontsize{9}{11}\selectfont,
  node distance=5mm,
  card/.style={
    rounded corners=4pt,
    align=flush left,
    text width=3.35cm,
    minimum height=3.3cm,
    inner sep=6pt,
    line width=0.8pt
  },
  stage/.style={card,draw=SumBlue,fill=SumBlue,text=white},
  gate/.style={card,draw=PosRed,fill=PosRed,text=white},
  arr/.style={
    -{Latex[length=2.2mm,width=1.6mm]},
    draw=VSCLArrow,
    line width=1pt,
    shorten >=1pt,
    shorten <=1pt
  },
  review/.style={arr,draw=PosRed,dashed},
  repair/.style={arr,draw=VSCLMuted,dashed},
  edge label/.style={
    font={\sffamily\fontsize{8}{10}\selectfont
          \hyphenpenalty=10000\relax},
    text=VSCLMuted,
    align=center,
    inner sep=2pt
  },
  band/.style={
    rounded corners=4pt,
    fill=VSCLInk,
    text=white,
    align=flush left,
    inner xsep=6pt,
    inner ysep=6pt
  }
]

\node[stage] (src) {
  \textbf{1. List the chapter's claims}\\[5pt]
  Agents list what must be formalized. Other agents check the chapter
  for anything the list missed.
};
\node[stage,right=of src] (packet) {
  \textbf{2. Read each claim carefully}\\[5pt]
  Record what it asserts and every condition it needs, including those
  stated only in the proof.
};
\node[stage,right=of packet] (target) {
  \textbf{3. Write the claim in Lean}\\[5pt]
  Write a precise statement and the definitions it uses. Leave proving
  it for later.
};
\node[gate,right=of target] (audit) {
  \textbf{4. Check the translation}\\[5pt]
  Other agents check whether the Lean statement and its definitions
  express the book's claim.
};

\node[stage,below=12mm of audit] (proof) {
  \textbf{5. Prove the accepted statement}\\[5pt]
  Use earlier results or write a new proof. Keep the accepted statement
  and its definitions unchanged.
};
\node[gate,left=of proof] (check) {
  \textbf{6. Check the proof}\\[5pt]
  Lean checks the proof. Software also checks for changed statements
  or definitions, missing steps, and extra axioms.
};
\node[gate,left=of check] (chapter) {
  \textbf{7. Review the whole chapter}\\[5pt]
  Check every claim on the original list, organize the results, and
  rebuild the library.
};
\node[stage,left=of chapter] (lib) {
  \textbf{8. Add the results to VSCL}\\[5pt]
  Make the proved results available to later proofs, with links back
  to the book.
};

\draw[arr] (src) -- (packet);
\draw[arr] (packet) -- (target);
\draw[arr] (target) -- (audit);
\draw[arr] (audit) --
  node[edge label,left=3pt] {Accepted} (proof);
\draw[arr] (proof) -- (check);
\draw[arr] (check) -- (chapter);
\draw[arr] (chapter) -- (lib);

\draw[review]
  (audit.north) -- ++(0,6mm)
  -| (target.north);
\node[edge label,above=1pt]
  at ($(target.north)!0.5!(audit.north)+(0,6mm)$)
  {Wrong translation? Revise it and review again.};

\draw[review] (check.north) -- (target.south);
\node[edge label,left=3pt,text width=5.4cm]
  at ($(check.north)!0.5!(target.south)$)
  {Statement or definition wrong?\\
   Revise it and repeat the review.};

\draw[repair]
  (check.south) -- ++(0,-6mm)
  -| (proof.south);
\node[edge label,below=1pt]
  at ($(check.south)!0.5!(proof.south)+(0,-6mm)$)
  {Proof fails? Repair the proof, not the statement.};

\path let \p1=(src.north west), \p2=(audit.north east) in
  node[band,anchor=south west,
       text width={\x2-\x1-12pt},
       font=\sffamily\bfseries\fontsize{10}{12}\selectfont]
  at ($(src.north west)+(0,12mm)$)
  {How the skill formalizes a chapter};

\path let \p1=(lib.south west), \p2=(proof.south east) in
  node[band,anchor=north west,text width={\x2-\x1-12pt}]
  at ($(lib.south west)+(0,-12mm)$)
  {A coordinator assigns tasks and keeps track of unfinished work.
   Agents search VSCL and Mathlib for useful results before writing
   new statements or attempting proofs.};

\end{tikzpicture}%
}
\caption{The skill's workflow for new material.}
\label{fig:skill}
\end{figure}

We use an internal formalization skill to build \vscl{}. Figure~\ref{fig:skill} shows the workflow for new material. The skill itself is not part of the public
release.

\textit{Listing the chapter's claims.} Agents read a
fixed copy of the chapter and list its mathematical claims. If the task covers
only a selected part of a chapter, that scope must be stated at the outset.
Other agents check the list against the source for omissions. Once agreed,
the list fixes the work to be considered for formalization. An agent cannot drop a difficult
claim and report the chapter as complete.
Unresolved claims remain on the list and count as unfinished work.

\textit{Recording each claim.} For each claim, the agents record the relevant
passages, the mathematical objects, the assumptions, and exactly what is
asserted. This record, called a \emph{claim packet}, includes conditions
stated only in the source's proof or surrounding text.

\textit{Reusing earlier results.} The agents search \vscl{} and a version of Mathlib before writing the Lean statement or proving
the claim or any results needed for its proof. An existing result can be
reused only when its definitions, assumptions, constants, and conclusions
supply what the proof needs.

\textit{Writing and reviewing the statement.} An agent writes the claim as a
Lean statement, together with the definitions it uses, but it does not yet
attempt a proof. Agents other than the author compare the statement with
the source, checking separately whether the Lean statement implies the
source claim and whether the source claim implies the Lean statement.
Our review adapts the semantic-correctness method of Meek et
al.~\cite[Section~3.2.1]{openmath2026}. In the more extensive reviews, one
agent also translates the Lean statement back into mathematical prose
without seeing the source, and another compares that translation with the
source. The review includes the definitions used by the statement,
including those imported from other files, so that a correct-looking
formula cannot conceal a different mathematical object or algorithm.
The current selective procedure uses one, two, or additional separate
reviews according to the claim's risk. Only a statement that passes review is sent to a proving agent.
These implication checks are judgments by the reviewing agents, not formal
proofs of equivalence, separate reviews can still share the same mistake.

\textit{Proving and checking the statement.} The proving agent uses existing
theorems or writes a new proof. Automated checks require it to prove the
reviewed statement, without changing the definitions on which that
statement depends, and using only the permitted axioms. Accepted proofs
may contain neither placeholders for missing steps nor new axioms. A
failed proof is repaired with the statement and definitions unchanged.
If a statement or definition must be corrected, it must be reviewed again
before proving continues. A change to a shared definition also requires
reconsidering the statement reviews of claims that depend on its meaning, i.e.,
rebuilding their proofs is not enough. Lean checks the proof of the formal
statement.

\textit{Checking coverage and publishing results.} A coordinator assigns tasks, tracks which
results depend on which others, and records what has been completed and
what remains. Before reporting completion of an audit, the workflow checks
every claim in the agreed scope, its statement reviews and proof status,
the placement of accepted results in the library, and a full build.

\textit{Incorrect claims and unfinished proofs.} If a printed claim is false,
the workflow requires a formal counterexample and a separately named,
proved correction. The counterexample and correction are recorded
separately, and neither counts as a proof of the original statement.
A statement with extra assumptions, a narrower domain, or a weaker
conclusion is reported with that change.

\subsection[Building on VSCL]{Building on the library}
\label{sec:randnla}

This subsection gives evidence that \vscl{} extends from textbooks to research results. A
reader may fairly ask whether a large library written by AI agents is worth anything. We answer
with a case study in which we prove results of a research survey with the library and publish the
statements so that a reader can check them against the survey.

\paragraph{A case study in RandNLA.}\label{par:randnla-case-study}
Based on the survey of Drineas and Mahoney~\cite{drineasmahoney2016cacm}, the case study adds to
the library random sampling for matrix approximation and least-squares problems, two core problems
in randomized numerical linear algebra (RandNLA)~\cite{murray2023randnla}.
It covers sampling individual
matrix entries, sampling and rescaling rows, and preprocessing with
random projections (Table~\ref{fig:coverage-randnla}). The associated
analysis includes spectral and Gram-matrix approximation,
leverage-score sampling and subspace preservation, sketched
least-squares objectives, and low-rank approximation. 
Table~\ref{tab:randnla-use} provides information about how the formalized results for this case study depend on the rest of the library. In particular, each proof
depends on 91 to 589 library theorems. The case study's record on our website lists each with its place
in the survey, its Lean text, and links to its definitions and proof~\cite{vsclrandnlarecord2026},
so a reader can compare it with the survey without reading a proof. 
\begingroup
\let\vsclOriginalCaption\caption
\renewcommand{\caption}{%
  \par\vspace{3pt}%
  \vsclOriginalCaption}
\begin{table}[!htb]
\centering
\begingroup
\normalfont\small
\setlength{\parindent}{0pt}
\setlength{\parskip}{0pt}
\newcommand{\covrow}[2]{%
  \noindent\begin{minipage}[t]{0.29\linewidth}\raggedright
  #1\strut\par
  \end{minipage}\hfill
  \begin{minipage}[t]{0.68\linewidth}\raggedright
  #2\strut\par
  \end{minipage}\par\vspace{4pt}}
\begin{minipage}{\linewidth}
\textbf{Drineas--Mahoney~\cite{drineasmahoney2016cacm}}\enspace\emph{RandNLA: Randomized Numerical Linear Algebra}\par
\vspace{4pt}\hrule height 0.6pt\vspace{5pt}
\covrow{Algorithm 1; equation (2)}{\emph{Entry sampling.} Rescaled entries and matrix approximation bounds.}
\covrow{Algorithm 2; equations (4)--(5)}{\emph{Row sampling.} Squared-norm probabilities, rescaling, and Gram-matrix approximation.}
\covrow{Equations (6)--(7)}{\emph{Subspace preservation.} Leverage scores and concentration bounds for changes in vector lengths.}
\covrow{Equation (8)}{\emph{Sketched least squares.} Residual objectives, approximation guarantees, probability bounds, and rounding error.}
\covrow{Equation (9)}{\emph{Low-rank approximation.} Column sketches, rank factorizations, projectors, and approximation conditions.}
\covrow{Algorithm 3}{\emph{Random projections.} Left, right, and two-sided transforms; exact products and rounding error.}
\hrule height 0.4pt
\end{minipage}
\endgroup
\caption{Algorithms and equations from~\cite{drineasmahoney2016cacm} on sampling, random
projections, matrix approximation, and least squares.
Probability and approximation guarantees are developed apart
from rounding-error bounds, each under its own sampling,
conditioning, and arithmetic assumptions~\cite{lcm2026repo}.}
\label{fig:coverage-randnla}
\end{table}

\endgroup

\begin{table}[!htb]
\centering
\small
\begin{tabular}{@{}llr@{}}
\toprule
Result & In the survey & Library theorems \\
\midrule
Spectral-norm bound for element-wise sparsification & Eq.~(2) & 430 \\
Structural bound for low-rank approximation & Eq.~(9) & 184 \\
Subspace embedding by the randomized Hadamard transform & Alg.~3 & 565 \\
Least squares with the randomized Hadamard transform & p.~87 & 589 \\
Least squares with approximate leverage scores & p.~87 & 451 \\
Sparse subspace embedding & Alg.~3 & 91 \\
\midrule
All 35 published proofs & & 1{,}053 \\
\bottomrule
\end{tabular}
\caption{Six results of the RandNLA survey~\cite{drineasmahoney2016cacm} with proofs in \vscl{}~\cite{lcm2026repo}, and the number of theorems of the library that
each proof depends on, directly or through a chain of references in Lean's elaborated statements
and proofs, as in the audit behind Table~\ref{tab:library-reuse}. Theorems of Mathlib and the 57
theorems in the case study's own files are not counted. The proofs share many theorems, so the rows
add up to more than the last. For the entry-sampling result, we prove the bound for square matrices
using the thresholded algorithm of Drineas and Zouzias. The survey
omits the thresholding step needed for this bound. For the Gram-matrix
result, we make the expectation explicit, as specified in the survey's
prose, and define the sampler to return a zero sketch when the input
matrix is zero.}
\label{tab:randnla-use}
\end{table}

\section[\veritas{}]{\veritas{}}
\label{sec:veritas}

Our audit (Section~\ref{sec:audit}) called for an assessment of NLA kernels that does not depend on a fixed test set. 
The error bounds formalized in \vscl{} (Section~\ref{sec:vscl}) specify the largest error a kernel may produce on every input they cover.
Our \veritas{} project is designed to assess any NLA kernel against its specification, which holds the error bounds that apply to the kernel, the required treatment of exceptional values, and the tolerances that the tests use. 
\veritas{} is designed to take each part from \vscl{} once \vscl{} has formalized it.
Today, \vscl{} proves its error bounds in the standard model of floating-point arithmetic~\cite[Section~2.2]{higham2002asna}, in which each operation returns the exact result up to a relative error of at most the unit roundoff and no operation overflows or underflows. 
Such a bound applies to a run of a kernel only where the kernel's arithmetic meets these assumptions, and establishing when it does is future work (Section~\ref{sec:bridge}).

\begin{figure}[!htb]
\centering
\begin{tikzpicture}[
  font={\sffamily\fontsize{8.5}{10.5}\selectfont\hyphenpenalty=10000\relax},
  card/.style={rounded corners=4pt,align=flush left,text width=2.4cm,minimum height=1.7cm,
               inner sep=6pt,line width=0.8pt},
  stage/.style={card,draw=SumBlue,fill=SumBlue,text=white,minimum height=2.25cm},
  source/.style={card,draw=black!40,fill=black!4,text=black!85,text width=2.0cm},
  library/.style={card,draw=SumBlue,fill=white,text=black!85,align=center,
                  text width=4.2cm,minimum height=1.0cm},
  arr/.style={-{Latex[length=2.3mm,width=1.7mm]},draw=black!60,line width=1pt,
              shorten >=1.5pt,shorten <=1.5pt},
  lab/.style={font=\sffamily\fontsize{8}{9.5}\selectfont,text=black!85,align=center,inner sep=3pt}]
\node[source] (kernel) at (0,0) {\textbf{Kernel}\\[2pt] human- or AI-written};
\node[stage] (s1) at (3.22,0) {{\fontsize{13}{13}\selectfont\bfseries 1}\enspace\textbf{Recognize}\\[2pt] what the kernel computes and how};
\node[stage] (s2) at (6.64,0) {{\fontsize{13}{13}\selectfont\bfseries 2}\enspace\textbf{Test}\\[2pt] against error bounds, exception rules and tolerances};
\node[stage] (s3) at (10.06,0) {{\fontsize{13}{13}\selectfont\bfseries 3}\enspace\textbf{Certify}\\[2pt] record what was checked and found};
\node[library] (vscl) at (6.64,-3.05) {\textbf{\vscl{}}\\ machine-checked error bounds\\ {\fontsize{7.5}{9}\selectfont exceptional values planned}};
\begin{scope}[on background layer]
  \node[draw=SumBlue!45,fill=SumBlue!7,rounded corners=7pt,line width=0.6pt,inner sep=7pt,
        fit=(s1)(s2)(s3)] (band) {};
\end{scope}
\node[anchor=south west,font=\sffamily\bfseries\fontsize{9}{10}\selectfont,text=SumBlue]
  at ([xshift=3pt,yshift=1pt]band.north west) {\veritas{}};
\path[draw=PosRed,fill=white,line width=0.9pt,line join=round]
  (12.07,-1.1) -- (14.87,-1.1) -- (14.87,0.7) -- (14.47,1.1) -- (12.07,1.1) -- cycle;
\path[draw=PosRed,fill=PosRed!18,line width=0.9pt,line join=round]
  (14.47,1.1) -- (14.47,0.7) -- (14.87,0.7) -- cycle;
\node[align=flush left,text width=2.3cm,inner sep=0pt,text=black!85] at (13.42,-0.05)
  {\textbf{\textcolor{PosRed}{Certificate\\ of trust}}\\[3pt] {\fontsize{7.5}{9}\selectfont checks, inputs and environment, outcomes}};
\coordinate (certw) at (12.07,0);
\draw[arr] (kernel) -- (s1);
\draw[arr] (s1) -- (s2);
\draw[arr] (s2) -- (s3);
\draw[arr] (s3.east) -- (certw);
\draw[arr] ([xshift=-5mm]vscl.north) -- node[lab,left] {error\\ bounds} ([xshift=-5mm]s2.south);
\draw[arr,dashed] ([xshift=5mm]s2.south) -- node[lab,right] {missing\\ analyses} ([xshift=5mm]vscl.north);
\end{tikzpicture}
\caption{The three steps by which \veritas{} is designed to assess a kernel, and the exchange between \veritas{} and \vscl{}.}
\label{fig:veritas-goal}
\end{figure}
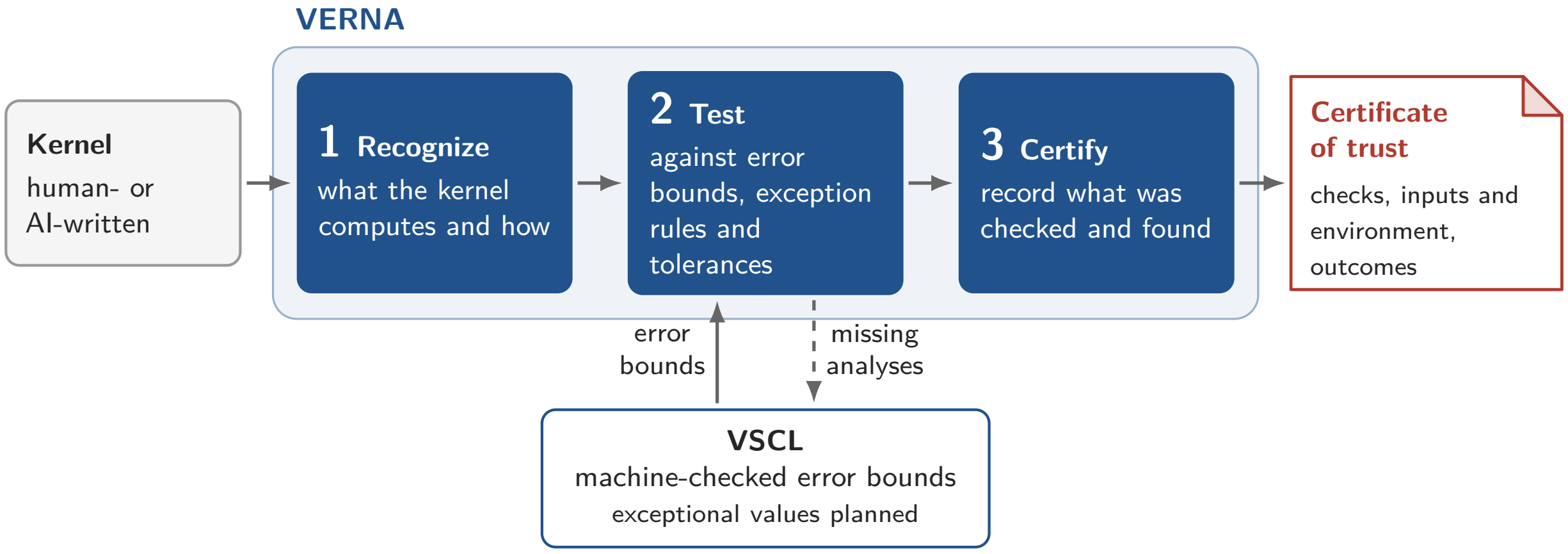

\veritas{} is designed to assess a kernel in three steps (see Figure~\ref{fig:veritas-goal}).
It first recognizes the kernel, meaning the operation it computes and the algorithm and arithmetic it uses.
It then tests the kernel against its specification, on inputs where unstable methods tend to fail, such as ill-conditioned or nearly rank-deficient matrices. The error bounds in the specification depend on the input, e.g., on its order and norm, so the test evaluates each of them on the input at hand. Exceeding a bound proved in \vscl{} is a reason to reject the kernel, but meeting it is not enough to accept the kernel, because worst-case bounds are often loose. Until \vscl{} covers exceptional values, they are checked as Demmel and co-authors propose~\cite[Section~9]{demmelgrading2026}.
Finally, \veritas{}
issues a certificate of trust, which records the checks performed, the inputs and the execution environment, and the outcomes.
When a failed test is not the kernel's fault and no formalized bound accounts for it, the case returns to \vscl{} as a task to formalize the missing analysis for that kernel.

Our audit carried out a simplified form of the testing step in \veritas{}.
It checked the competition's requirements, finite output and the two tolerances of Table~\ref{tab:qr-contract}, on inputs from the competition's generator.
Its record identifies the checkers, the inputs, and the environments, and gives the outcome of every recorded evaluation, as a certificate would~\cite{melnichenko2026auditrecord}.

The goal of \veritas{} is to make assessment against a specification routine for any NLA kernel, no matter who or what wrote it.
As Section~\ref{sec:audit} notes, following Krentsel and co-authors~\cite{krentsel2026reality}, in NLA the largest errors are fixed by the error bounds of numerical analysis, the arithmetic and its exceptional values are specified, and reference implementations exist. We believe that this makes NLA a domain where software verification is plausible and must be pursued. \vscl{} already extends beyond NLA into other areas of numerical analysis, and as it grows we aim to assess numerical kernels more broadly. We see this work as a precursor to verifying software more generally.
This paper states the goal of \veritas{} and leaves its design to a separate publication.

\section{What lies ahead}
\label{sec:conclusion}

As \vscl{} grows, our aim is to develop reusable stability analyses that can be adapted to changes in algorithms, precision, and implementation.
In this section, we describe what lies ahead. 
We first describe priorities for formalization, from error bounds
for basic NLA kernels to guarantees for complete algorithms.
We then discuss our aim to connect these results to concrete arithmetic
and compiled implementations, including the treatment of underflow
and exceptional values.
We then consider open-weight models as a way to reduce the cost
of formalization.
We conclude by returning to the question of the title and setting out six design goals inspired by Turing's work; these are principles that should guide efforts at developing, reviewing, and using formalized results to assess numerical software.

\subsection{Formalization targets in stability analysis}
\label{sec:analysis-value}

We aim to bring the broader vision of AI-assisted proof
adaptation~\cite[Section~4]{tao2025machine} to numerical stability analysis.
Building on \vscl{}, we will formalize stability proofs that can be adapted to changes in implementation or precision. This will help establish sharp guarantees for individual implementations of numerical kernels and complete algorithms at scale. 

Numerical analysts already adapt earlier arguments to new algorithmic
and arithmetic settings, see e.g.~\cite{connolly2021stochastic,melnichenko2025cqrrpt,
connollyhigham2023qr}. \vscl{} will provide a basis for making this practice systematic by organizing local error estimates and propagation arguments as reusable formal components.
The library reuse analysis in Appendix~\ref{app:benchmark_details}
supports this perspective.

Our choice of future formalization targets will therefore be guided
by the sharpness and applicability of their guarantees and the
adaptability of their proofs. These targets include analyses that
derive algorithm-level error bounds from accuracy guarantees for
constituent operations~\cite{stabilityS002,carson2026trackable,
afnButtari2026GMRES}, as well as analyses of specific algorithms whose
proofs separate local error estimates from their
propagation~\cite{yang2021mixedqr,stabilityS032,carson2022blockgs}.
Three complementary formalization priorities follow.

We will begin by formalizing stability analyses of basic NLA kernels,
whose error estimates support many larger algorithms. We will prioritize
proofs that exploit properties of rounding and the operands to obtain
sharper bounds or less restrictive applicability conditions. Refined
rounding-error analyses provide natural starting
points~\cite{jeannerodrump2013,rump2014factorizations,carson2026trackable}. Such bounds can depend on implementation choices, including accumulation
schemes and precision assignments~\cite{carsonma2026ozaki,blanchard2020blockfma,stabilityS072}.
Further refinements can exploit matrix structure, restrictions on entry
magnitudes, or probabilistic models of inputs and rounding
errors~\cite{afnHighamMary2020,huang2024average}.
Such refinements can be essential for obtaining informative guarantees, particularly in low precision floating-point formats.

For applications in numerical PDEs, optimization, and machine learning, we will formalize results that derive conditioning and spectral properties from the underlying problem, its discretization, and preceding computations~\cite{mardalwinther2011preconditioning,afnBalabanovGrigori2022,kamenski2014conditioning,pazner2023loworder,louart2018random,fanwang2020spectra}.
 We will also formalize results on algorithmic behavior. For iterative
methods for linear systems and eigenvalue problems, these include projection and orthogonality relations,
optimality, and convergence~\cite{saad2003iterative,saad2011eigenvalue,afnGraham2017}.
For nonlinear equations, we will study how inexact computations
affect error bounds of second-order methods~\cite{tisseur2001newton,stabilityS080}.

Building on these results, our third priority is to formalize how
component errors jointly affect convergence and final
accuracy~\cite{stabilityS002,stabilityS007,afnButtari2026GMRES}.
The aim is to identify sufficient accuracy requirements for individual
components, accounting for how errors introduced at different stages
are amplified, accumulated, or reduced by the
algorithm~\cite{carsonhigham2018,afnMcCormick2021,
afnElArar2025,afnBeuzeville2026,afnBalabanov2026QlessQR}.
This includes assessing how changes in preconditioning, update formulas,
or precision assignments alter these requirements and the resulting
guarantees~\cite{afnMeier2024,stabilityS027,ketcheson2014internal,stabilityS082}.

\subsection{From numerical methods to mathematical software}
\label{sec:bridge}

Fine-grained numerical analysis becomes increasingly challenging as algorithms grow in complexity to exploit specialized hardware~\cite{uchino25,haidar2020,flashattention3, fu2024}. At the same time, hardware–software co-design drives the development of more specialized architectures, whose internal arithmetic behavior may not be fully specified or exposed to software~\cite[Sections~9.7.14.5 and 9.7.15.4]{nvidiaPTX}~\cite{valpey2025}. As both trends seem to accelerate with the adoption of AI-assisted design, incorporating implementation-specific choices into algorithm specifications is needed for guarantees that reflect actual software behavior. 

For an IEEE binary floating-point format \(F\), let \(a\circ_F b\) denote an elementary arithmetic operation on operands \(a,b\in F\). Under appropriate rounding, underflow, and overflow assumptions, its total error can be bounded in terms of arithmetic rounding  and input representation errors~\cite[Sections 2.3.1 and 2.3.2]{muller2018}. Constructing similar guarantees for numerical kernels that rely on arbitrary combinations of storage and arithmetic schemes requires an automated approach to formalizing their computational behavior. This involves generating and proving composition theorems that propagate error bounds across operations.

The LAProof~\cite{laproof2023} and VCFloat2~\cite{vcfloat2_2024} libraries, described in Appendix~\ref{sec:algo-level}, automate the derivation of roundoff error bounds for floating-point expressions and include support for nonstandard formats and external functions. Valpey and co-authors~\cite{valpey2025} formalize the arithmetic behavior of NVIDIA Tensor Cores,  using SMT reasoning to investigate accumulation precision, rounding, and the numerical behavior of mixed-precision matrix multiplication algorithms. These efforts primarily rely on manually developed formalizations and established automated reasoning techniques. LLM-assisted systems  demonstrate progress in automating formalization \cite{kan2026_aria} and theorem-proving  \cite{hilbert_formal_proofs, goedel_prover_v2, kozyrev2026}
 but do not specifically target the construction and certification of error bounds. Currently, \vscl{} proves algorithm-level bounds in the standard model which does not account for overflow or underflow.
By utilizing libraries that formalize IEEE and non-standard arithmetic semantics \cite{floatlib2026, tunnell2026ieee}, \vscl{} can leverage its collection of formalized theorems to support specification of numerical kernels, and verification of correctness theorems over a wide range of storage and arithmetic types. VERNA complements this process by identifying gaps in existing formalizations and guides proof generation. It does so by adding results to a kernel's specification as \vscl{} formalizes them, while identifying the cases that are not covered by any formalized bound.

Applying an analysis theorem to an execution requires the execution to meet the theorem's finiteness, domain, underflow, and no-overflow or no-saturation conditions. %
Compiler transformations such as instruction reordering, arithmetic re-association or operation contraction enabled by underlying hardware also remain explicit assumptions that the certificate of trust can simply record. Furthermore, while IEEE 754 specifies the behavior of exceptional values, the standard relative error model does not generally account for exception propagation. Extending numerical analysis theorems to describe the propagation of exceptional values, including NaNs and infinities, through numerical algorithms remains an important future goal. %

Formal analysis of numerical algorithms with hardware-aware specifications does not, by itself, establish software correctness. Bridging this gap requires a formal correspondence between algorithms and their compiled implementations. Existing efforts demonstrate this connection; LAProof, VCFloat and VCFloat2 connect floating-point error analysis, with C program verification. 
A complementary approach starts from the compiled program instead, to reconstruct its arithmetic operations and data dependencies. Prior work in~\cite{lee16, lee18} has demonstrated formal error analysis of compiled floating-point implementations, particularly optimized mathematical functions, while dynamic tools such as Herbgrind~\cite{herbgrind} and Verrou~\cite{verrou} track floating-point errors-propagation in compiled programs and investigate the effects of rounding perturbations respectively. By incorporating elements from these approaches, \vscl{} can extend numerical formalization toward software verification, connecting numerical theorems to implementations by specifying  the actual low-level transformations and hardware behavior. 

\subsection{Open-weight models}
\label{sec:open-weights-summary}

Formalizing material from a single book requires months of sustained
model use. Our six subscriptions at US\$200 per month cost US\$1,200, and personal subscriptions
impose usage limits~\cite{openai2026codexpricing}. API access supports
larger workloads, but one of our 13-hour runs cost approximately
US\$1,300. Sustaining this spending rate for one month would cost
approximately US\$70,000. Subscriptions cost far less than sustained API use, so access at subscription prices may not last, and a project that depends on it could face costs closer to the API rate. These expenses motivate methods that reduce the cost of extending \vscl{}.

Open-weight models offer one possible approach. We are pursuing two
directions. The first is to run the formalization skill described in
Section~\ref{sec:skill} using a self-hosted, open-weight coding agent.
The second is an orchestrator--prover workflow that separates statement
construction and review from proof generation, so that two smaller,
specialized models can share the work. As a first experiment, DeepSeek-V4.1-Flash-1017 and GPT-6.1 Sol each attempted to formalize the Appetizer and Chapter 1 of \emph{High-Dimensional
Probability}~\cite{vershynin2018hdp} from the source PDF
(Table~\ref{tab:standalone-hdp}). Both produced substantial proofs
that Lean accepted. Appendix~\ref{sec:open-weights}
gives the details.

\begin{table}[!htbp]
\centering
\small
\renewcommand{\arraystretch}{1.25}
\setlength{\tabcolsep}{6pt}
\begin{tabular*}{\textwidth}{@{\extracolsep{\fill}}llccccc@{}}
\hline
\noalign{\vskip 5pt}
Model & Material & \shortstack{Active time\\(h:mm:ss)}
& \shortstack{Input\\(M)}
& \shortstack{Cached input\\(M)}
& \shortstack{Output\\(M)}
& \shortstack{Cost\\(US\$)} \\
\noalign{\vskip 5pt}
\hline
DeepSeek-V4.1-Flash-1017 & Appetizer & 0:48:35 & 81.426 & 79.156 & 0.215 & 1.06 \\
DeepSeek-V4.1-Flash-1017 & Chapter~1 & 3:01:40 & 310.967 & 301.870 & 0.446 & 3.18 \\
GPT-6.1 Sol & Appetizer & 0:39:22 & 45.843 & 44.380 & 0.285 & 10.21 \\
GPT-6.1 Sol & Chapter~1 & 0:52:34 & 83.944 & 82.256 & 0.374 & 16.47 \\
\hline
\end{tabular*}
\caption{Standalone formalization experiments without access to \vscl{} or prior
formalizations of the material. Token counts are reported in millions (M).
Cached input is included in total input. DeepSeek-V4.1-Flash-1017 costs are
recorded API charges. GPT-6.1 Sol ran through ChatGPT
sign-in~\cite{openai2026codexpricing}, so its costs are estimates calculated
from recorded tokens, observed caching, and the applicable API
rates~\cite{openai2026apipricing}. The code for these experiments is available
in~\cite{openweight2026repo}.}
\label{tab:standalone-hdp}
\end{table}

Table~\ref{tab:standalone-hdp} shows that DeepSeek-V4.1-Flash-1017 produced substantial
Lean proofs at low API cost, US\$4.24 across the two units, compared
with an estimated API-equivalent cost of US\$26.68 for GPT-6.1 Sol. GPT-6.1 Sol required less active time, approximately 1~hour 32~minutes
in total versus 3~hours 50~minutes for DeepSeek-V4.1-Flash-1017, and consumed fewer
input tokens in both runs. Approximately 97--98\% of input tokens
were cached in all four runs. These figures must be interpreted
alongside the differences in workflows, operator assistance, and
final coverage described in Appendix~\ref{sec:open-weights}.
For example, DeepSeek-V4.1-Flash-1017 left the covering-sharpness claim unproved,
whereas GPT-6.1 Sol proved a precise statement supplied after clarification.
The reported costs also exclude operator preparation, specification
research, and final audits.

These first runs suggest that open-weight models could carry part of the formalization work, and we intend to pursue both directions.

\subsection{What would Turing do?}
\label{sec:wwtd-principles}

The audit of the competition in Section~\ref{sec:audit} illustrates questions that arise whenever a numerical method is turned into software: what accuracy the method can deliver, which assumptions that accuracy requires, and how to establish that a computation meets them. These questions were already central to the early development
of digital computing.
Von Neumann and Goldstine published a rounding-error analysis of matrix inversion in 1947~\cite{vonneumann1947inverting}.
Turing's 1948 paper analyzed rounding errors in methods for solving linear systems and inverting matrices~\cite{turing1948rounding}.
In 1949, Turing described a method for checking programs using assertions about program states and a separate termination argument~\cite{turing1949checking}. 
We draw on these and other parts of his work to formulate six design goals for our research program. 
Many are established ideas in numerical analysis and formal verification. 
Their role here is to guide which results we formalize, how we review them, and how we use them to assess numerical software.

\paragraph{1. Review whether the formal statement expresses the intended claim.}
In Section~9 of \emph{On Computable Numbers}, Turing supported
the claim that his formal definition captured the intuitive
notion of computability by analyzing how a person carries out
a calculation on paper~\cite{turing1936computable}.
Turing was asking whether his formal definition captured
everything that could be computed.
Our concern is whether a formal statement faithfully represents
the mathematical claim expressed by its informal source.
These tasks differ in scope, but both require assessing the
correspondence between a formal description and its intended
informal meaning.
In autoformalization, a checked proof establishes a formal
statement under its assumptions and axioms, but it does not by
itself establish that the statement faithfully represents its
informal source~\cite{openmath2026,choirpr143}.
To assess this correspondence in our setting, human guidance
and targeted expert checks could complement the existing agent
reviews by comparing the formal and source statements. The RandNLA case study publishes its statements apart from their proofs for exactly such checks (Section~\ref{par:randnla-case-study}).

\paragraph{2. State the representation and arithmetic assumptions.}
In his analysis of Jordan's method, Turing imposed bounds on
errors in individual computed quantities and noted that the
precision needed to meet those bounds depends on the values
involved~\cite[Section~11]{turing1948rounding}.
Turing was deriving error estimates under specified
arithmetic assumptions.
Our task is to formalize such estimates and establish when
the arithmetic of a numerical kernel satisfies their
assumptions.
In both cases, applying an error estimate to a computation
requires justifying the bounds assumed for its individual
operations.
For example, the usual relative-error bound does not cover
a nonzero exact result that rounds to zero.
To apply an analysis using that bound, we must establish that
this case cannot occur on the covered inputs.
To cover such cases, the analysis must account for their
errors.
We must therefore specify how stored values relate to the
exact quantities in the theorem, the precision and rounding
rule used by each operation, and the treatment of underflow,
overflow, and exceptional
values~\cite{ieee754_2019,posit2022,p3109}.
Using further properties of round-to-nearest arithmetic can
yield sharper bounds for inner products and matrix
factorizations~\cite{jeannerodrump2013,rump2014factorizations}.
Connecting \vscl{}'s abstract rounding model to concrete
arithmetic remains future work (Section~\ref{sec:bridge}).

\paragraph{3. Distinguish sensitivity, approximation error, and rounding error.}
Turing studied how perturbations in the coefficients of a
linear system affect its solution.
He also discussed iterative refinement, which uses the
residual of an approximate solution to compute a correction.
His discussion of rounding errors emphasizes the order of
operations and the precision of intermediate
results~\cite[Sections~8--10]{turing1948rounding}.
For our project, this distinction should guide how we state
and use formalized numerical results.
For a nonsingular linear system, conditioning describes how
perturbations in the data affect the exact solution.
An approximation bound can describe the error remaining
after finitely many iterations even in exact arithmetic.
A rounding-error bound describes the effects of the
specified rounded operations~\cite{higham2002asna}.
A numerical accuracy guarantee may require combining these
results.
The formalization should make explicit which errors are
bounded, how problem sensitivity enters those bounds, and
which assumptions are needed.

\paragraph{4. Build the argument from smaller claims with explicit assumptions.}
In \emph{Checking a large routine}, Turing proposed checking
a routine through assertions that could be verified
individually and used together to establish its correctness.
He used assertions about intermediate program states and
gave a separate termination
argument~\cite{turing1949checking}.
Turing connected claims about successive program states.
In our numerical analyses, we connect mathematical error
bounds for operations and intermediate results.
In both cases, the larger argument requires establishing
that the assumptions needed at each step follow from the
available facts.
For example, \vscl{}'s matrix-multiplication error proof uses
a matrix-vector bound, which in turn uses a dot-product bound.
The latter uses the summation analysis and the
addition-rounding model (Section~\ref{sec:library}).
These results can be reused without repeating their proofs,
but each application must establish their hypotheses.
This composition of mathematical results is already present
in \vscl{}.

\paragraph{5. Where possible, certify the particular answer returned.}
Turing distinguished assessing the accuracy of an answer
already computed from estimating the accuracy of a proposed
method. For a computed inverse \(B\), he used the residual \(I-AB\)
to bound its error~\cite[Section~10]{turing1948rounding}.
Our corresponding goal is to develop formally justified
checks that establish an accuracy guarantee for a particular
numerical output, without requiring a correctness proof
for the program that produced it.

For example, let \(A,B\in\mathbb{R}^{n\times n}\) be nonsingular,
let \(b\in\mathbb{R}^n\setminus\{0\}\), and let \(x=A^{-1}b\).
Using \(B\) as a left preconditioner, a checker can bound the
relative forward error of an approximation
\(\widehat{x}\in\mathbb{R}^n\) by
\[
\frac{\|x-\widehat{x}\|_2}{\|x\|_2}
\leq
\operatorname{cond}_2(BA)\,
\frac{\|B(b-A\widehat{x})\|_2}{\|Bb\|_2}.
\]
To certify that \(\widehat{x}\) meets the requested relative error
tolerance, the checker must establish an upper bound on
\(\operatorname{cond}_2(BA)\), an upper bound on
\(\|B(b-A\widehat{x})\|_2\), and a lower bound on
\(\|Bb\|_2\). These bounds and their final combination must account
for rounding errors in all numerical evaluations.

We aim to apply this principle to certify output accuracy and establish
algorithm-level conditions required by stability or convergence analyses,
including bounds on conditioning and loss of orthogonality, spectral
constraints, and bounds on local errors and their amplification.

These certificates can complement the certificate of trust in \veritas{}, which records the checks performed for a given input and computed result (Section~\ref{sec:veritas}).

\paragraph{6. Allow flexible proof search under fixed acceptance criteria.}
In his discussion of learning machines, Turing observed that
a machine can follow valid logical rules while choosing
deductions that do little to help it reach its
goal~\cite[Section~7]{turing1950computing}.
Our workflow applies this distinction by allowing agents
to try different proof strategies while keeping the
requirements for accepting a proof fixed.
Once a formal statement and its definitions have been
reviewed against the source, proving agents choose lemmas,
explore arguments, and construct candidate proofs.
An accepted proof must establish the reviewed statement,
with its definitions unchanged and under the permitted
axioms.
Lean checks whether the resulting proof establishes the
formal statement.

Together, these goals describe a program for developing reusable, checkable numerical guarantees as algorithms, representations, and precision choices change. The ambition is to combine numerical expertise with AI-assisted exploration and formal proof checking, so that faster generation of software is accompanied by better evidence of its reliability. \vscl{} is where such guarantees are proved, and \veritas{} is designed to assess NLA kernels against them.

\section*{Acknowledgments}

K.~Fountoulakis would like to acknowledge the support of the Natural Sciences and Engineering Research Council of Canada (NSERC). Cette recherche a \'et\'e financ\'ee par le Conseil de recherches en sciences naturelles et en g\'enie du Canada (CRSNG), [RGPIN-2019-04067, DGECR-2019-00147].

Sandia National Laboratories is a multimission laboratory managed and operated by National Technology \& Engineering Solutions of Sandia, LLC, a wholly owned subsidiary of Honeywell International Inc., for the U.S. Department of Energy’s National Nuclear Security Administration under contract DE-NA0003525.

This paper describes objective technical results and analysis. Any subjective views or opinions that might be expressed in the paper do not necessarily represent the views of the U.S. Department of Energy or the United States Government.

\subsection*{Use of AI tools}
Not surprisingly, given our interest in using modern AI tools for the future of numerical analysis and computational mathematics software, we did not shy away from their use.

The authors used Claude Code (Anthropic) \cite{anthropic2026claudecode} and Codex (OpenAI) \cite{openai2026codex} under human direction to gather and organize the literature, to draft and revise parts of the text, to check citations against their sources, to produce the figures from the data described in their captions, and to run the kernel audits of Section~\ref{sec:audit} and Appendix~\ref{sec:audit-details}. AI systems are also part of the work that the paper reports. \vscl{} is built with the workflow of Section~\ref{sec:skill}, in which AI agents write the Lean statements and proofs, other agents review the statements against their sources, and Lean checks every proof. In the reuse analysis of Appendix~\ref{app:benchmark_details}, a language model extracted and classed the proof tasks and a coding agent attempted them. Section~\ref{sec:open-weights-summary} and Appendix~\ref{sec:open-weights} report experiments with two further models.
The authors chose the questions and the claims, checked the numbers against their sources, and are responsible for all content.

\setlength{\bibsep}{0.6pt plus 0.3pt}
\renewcommand{\bibfont}{\scriptsize}
\begin{multicols}{2}
\bibliographystyle{abbrvnat}
\bibliography{refs}
\end{multicols}

\clearpage
\appendix
\section{Existing results and tools for numerical methods}
\label{sec:prior}

In the coming subsections, 
we first discuss work that establishes guarantees for algorithms and programs (Section~\ref{sec:algo-level}),
we then turn to their arithmetic foundations (Section~\ref{sec:arith}) and the AI-assisted workflows that can help build larger libraries
(Section~\ref{sec:ai-pipelines}).

\subsection{Numerical guarantees for algorithms and programs}
\label{sec:algo-level}

\paragraph{Theorems about whole algorithms.} 
Roux formally proved a bound for the floating-point Cholesky factor and used it in a check of positive definiteness \cite{roux2016}.
LAProof \cite{laproof2023}, a Coq library of accuracy proofs for linear algebra programs, proved error bounds for inner products, matrix-vector products, and matrix products using a rounding model that accounts for underflow.
It also connected its sparse matrix-vector product to C code through VCFloat \cite{vcfloat2_2024}, a Coq tool for floating-point reasoning, and the Verified Software Toolchain \cite{appel2011vst}.
The same group verified the Jacobi iteration in floating point, including its stopping test and an overflow check tied to the computed residual \cite[Sections~7--8]{tekriwal2023jacobi}.
Boldo and co-authors verified a finite-difference program for the wave equation \cite{boldo2013wave}.
Immler verified an executable ODE solver that chooses its time steps adaptively and computes enclosures accounting for both discretization and rounding errors \cite{immler2012ode,immler2018lorenz}.
Akbarpour and Tahar proved recurrences for accumulated roundoff error in fast Fourier transforms and showed that a register-transfer-level hardware design implements its fixed-point specification \cite{akbarpourtahar2004fft}. NLA-lib is an open-source Lean~4 library for numerical linear algebra,
including formalized results on randomized methods and conjugate
gradient in finite precision~\cite{nlalib2026}.

\paragraph{Automatic tools.} Some tools bound the rounding error of a given program automatically.
They handle straight-line code, a fixed sequence of arithmetic operations without loops, and code whose loops have a fixed bound.
Four such tools are Gappa, FPTaylor, Daisy, and PRECiSA \cite{gappa2010,fptaylor2019,daisy2018,precisa2024}. Gappa, from the group behind the Flocq arithmetic library \cite{flocq2011}, produces proofs that the Coq proof assistant checks.
FPTaylor bounds errors with symbolic Taylor expansions and supports subnormal numbers, uncertain inputs, transcendental expressions, and mixed precision within a fixed computation \cite{fptaylor2019}.
Its version 0.9.3 produces proof certificates checked in the HOL Light proof assistant.
PRECiSA produces proofs in the PVS proof assistant and bounds errors when rounding changes a branch decision; it also supports modular analysis of function calls and bounded list computations \cite[Sections~2--4]{precisa2024}.
Daisy analyzes mixed-precision assignments and searches for equivalent expressions with smaller error bounds.
Its certificates are checked by FloVer, whose checker is formally verified in Coq and HOL4 and whose HOL4 implementation is compiled using the verified CakeML toolchain \cite{flover2018}.

\subsection{Formal foundations for floating-point arithmetic}
\label{sec:arith}

The rounding-error results surveyed in Section~\ref{sec:algo-level} depend on assumptions about how individual arithmetic operations are rounded.
Formal arithmetic libraries help justify those assumptions for particular formats.
Some also prove bounds for sums and dot products, as described below.

Floating-point arithmetic has been formalized in several proof assistants, including Coq and Lean.
In Coq, Flocq provides a rounding theory over formats of any radix and precision, and most later Coq work builds on it \cite{flocq2011}. In Lean, three libraries were published in 2026.
\begin{itemize}\setlength{\itemsep}{0.5em}
\item
FLoPS formalizes the IEEE~P3109 family \cite{p3109} of narrow formats, including saturation and stochastic rounding, and proves properties of two short algorithms, FastTwoSum and ExtractScalar \cite{flops2026}. Its arithmetic is executable inside Lean and proved to refine the bit-level specification for rounding, saturation, and the four arithmetic operations, and it proves bias and root-mean-square bounds for sums under stochastic rounding \cite{flops2026repo}.
\item
Tunnell's FP library formalizes IEEE binary32 and binary64 arithmetic and proves error bounds for algorithms built on it: sequential and pairwise sums, dot products, and their compensated versions \cite{tunnell2026ieee}.
Its executable interval checkers also certify range conditions and accuracy bounds uniformly over input domains for supported sums and dot products.
\item
FloatLib covers a wide range of formats: IEEE binary and decimal formats, posits, P3109 \cite{p3109}, eight-bit and block-scaled formats, and user-defined formats, each modeled at the level of bits \cite{floatlib2026}. 
\end{itemize}

\subsection{AI-assisted formalization and statement review}
\label{sec:ai-pipelines}

\paragraph{Agentic autoformalization.}
RepoProver formalized a graduate textbook with agents in separate roles: sketchers write statements with placeholder proofs, provers fill them in, and separate mathematical and engineering reviewers check source faithfulness and code quality \cite[Section~2.4]{repoprover2026}.
Atlas reports applying the same approach to further online material from MIT OpenCourseWare \cite{rammal2026atlas}, and open platforms and open-source provers now let anyone's agents contribute \cite{chen2026prove2me,kripner2026openprover}.
Choir \cite{choir2026} runs a project as an open protocol: an orchestrator publishes statements as issues, workers' agents claim and prove them, and an automated gate rebuilds the project and rejects new axioms, undeclared placeholders, and changes to published statements.
It also accepts declared reductions, in which a worker splits a statement into smaller ones that are left open.
Townsend's catalogue of open problems in NLA formalizes those of its problems that already have informal solutions \cite{townsend2026openproblems}.
It fixes the hash of each Lean statement and has two reviewers other than the author approve it before any Lean proof body is written.
OpenMath~\cite{openmath2026} reports proof-complete Lean formalizations for 84 of the 175 statements in Butcher's \emph{Numerical Methods for Ordinary Differential Equations}~\cite{butcher2016ode}.

\paragraph{Faithfulness of statements.}
Choir \cite{choir2026} and Townsend's catalogue \cite{townsend2026openproblems} both use the Lean comparator \cite{leancomparator}, which checks that a proof establishes the statement as-written.
It does not check the statement's faithfulness to the source.
The Formal Conjectures benchmark warns, for the same reason, that its statements may not capture the original conjectures \cite{firsching2026formalconjectures}.
Two newer efforts address faithfulness directly.
Palomar, a registry of Lean proofs opened in August 2026, checks each proof mechanically and adds an AI review of whether the formal statement renders the informal one fairly, which its policy says is not expert review \cite{tao2026palomar}.
Tau Ceti, a Lean library built downstream of Mathlib \cite{mathlib4repo2026} by AI agents,  treats definitions formalized by AI as untrusted at first and builds confidence in them through adversarial review and through the results later built on them \cite{taucetiproject2026}.

\clearpage
\section{GPU kernel audit: setup and detailed results}
\label{sec:audit-details}

The audit selected ranks one through five from GPU MODE leaderboards 773 and 774 (QR rounds one and two), 775 (symmetric eigendecomposition), and 776 (Cholesky).
All twenty entries passed the scored cases from their own rounds, and eighteen failed on further inputs, in a step we call discovery (Table~\ref{tab:audit-overview}).

In this appendix, the vendor routine is NVIDIA's routine for the operation, which for QR is the routine that the baseline calls, and a target is an entry or a control.
For QR, discovery covers the audit inputs of Section~\ref{sec:audit}, called the common grid here, and the threshold cases described below.

The audit record contains the verdicts, input specifications, submission hashes, checker hashes, and environment records for all four rounds, and a script that replays a QR evaluation~\cite{melnichenko2026auditrecord}.
The organizers publish their checkers and input generators~\cite{gpumode2026referencekernels}, and their submission dataset contains the entries for leaderboards 774--776~\cite{gpumode2026dataset,kernelbot2025}.
Round-one QR source is absent from that dataset, so its recorded verdicts cannot currently be replayed from those public sources alone.

\begin{table}[!htbp]
\centering
\small
\caption{Audit of the five highest-ranked entries in each round.
All fractions count entries.
Scored replay uses the seven first-round QR inputs with their original checker and the twelve second-round inputs with the second-round checker.
QR discovery uses the second-round checker throughout.
Cholesky and eigendecomposition failures require the vendor to pass the same valid input.
Data: audit record~\cite{melnichenko2026auditrecord}.}
\label{tab:audit-overview}
\renewcommand{\arraystretch}{1.15}
\begin{tabular}{@{}lrrr@{}}
\toprule
Competition round & \makecell[br]{Failed own-round\\scored cases} & \makecell[br]{Failed in\\discovery} & \makecell[br]{Failure reproduced\\in reruns} \\
\midrule
QR, round 1 & 0/5 & 5/5 & Not run \\
QR, round 2 & 0/5 & 5/5 & Not run \\
Cholesky & 0/5 & 5/5 & 4/5 \\
Symmetric eigendecomposition & 0/5 & 3/5 & 3/5 \\
\bottomrule
\end{tabular}

\end{table}

\subsection{Environment, baseline time, and failure criteria}
\label{sec:audit-env}

All original audit evaluations ran on NVIDIA B200 GPUs, the competition's graded hardware model.
The organizers' runner pinned PyTorch 2.12.0, whose dependencies fix Triton 3.7.0.
Our QR audit used PyTorch 2.12.1 and Triton 3.7.1, one patch release later, with CUDA 13.0 runtime libraries.
It also used Python 3.12 where the runner used Python 3.13, and a newer CUDA compiler.
Every audited QR entry built and ran in that environment.
The Cholesky and eigendecomposition audit reconstructed the published runner recipe, including PyTorch 2.12.0.
The audit record preserves the actual versions and hashes for each execution.

The winner's leaderboard time is 1.22~ms.
The median of 48 submissions of the baseline in the organizers' dataset is 131~ms~\cite{gpumode2026dataset}, which gives the factor of 107 in Section~\ref{sec:audit}.
The audit record lists these 48 submissions and the factor for each of the five audited entries~\cite{melnichenko2026auditrecord}.
Participants' posts quote larger factors against a baseline time of 419~ms~\cite{lutz2026qr,sankalp2026autoresearch}, which only two of those 48 submissions approach.

We applied both QR checker versions to every QR output because the organizers tightened the checker between rounds.
Both QR checkers form $Q$ in float32 from the returned Householder vectors and coefficients.
Both evaluate the factorization and orthogonality residuals in float64.
A QR discovery failure means that the second-round checker rejected the entry's output.
The record labels an output as a fallback when it matches the vendor routine's output bit for bit, and no such output failed.
Upper-triangular and diagonal inputs can produce such matches without a fallback.
Among the second-round entries, only the fourth-ranked one matched outside these two families, on 311 inputs.
QR failures were not rerun in fresh containers.

A case specifies a matrix family, order, batch size, conditioning setting, and seed.
For QR, the audit ran each scored case once at its listed seed.
A leaderboard run checks one to fifty inputs per case at seeds derived from the listed one.

\subsection{QR inputs and controls}
\label{sec:audit-qr-inputs}

The QR sweep kept the seven scored shapes, because the entries choose their code path by shape.
It used the organizers' ten families: dense, upper-triangular, diagonal, rank-deficient, nearly rank-deficient, clustered, banded, near-collinear, row-scaled, and mixed.
For each of the seven scored shapes, five conditioning settings
and three seeds gave a common grid of 1,050 test cases per target.
This grid supplies the audit evaluations in Table~\ref{tab:qr-second-round} and Figure~\ref{fig:qr-families}.
It holds 945 distinct inputs, because the generator ignores the conditioning setting for clustered matrices and raises it to at least 4 for row-scaled ones.
Repeated inputs never changed a verdict.

For each target and scored shape, the audit also searched for a failing condition number with \texttt{randsvd} matrices~\cite[Section~28.3]{higham2002asna}.
The audit record holds one failing input at each threshold found, which adds at most five failures to a target's discovery count.

Each problem had two controls, the vendor routine and a degraded version that rounds each input to a lower precision before it calls the vendor routine.
The degraded QR control rounds to bfloat16, which keeps 8 of the 24 significant bits of single precision.
It tests whether the checker catches a loss of accuracy, and it is the second control in Table~\ref{tab:qr-second-round} and Figure~\ref{fig:qr-families}.
The checker rejected 619 of this control's 1{,}050 outputs on the common grid, including all 300 at orders 32 and 176.
The QR tolerances grow with $n$ by design~\cite{gpumode2026linalg}.
At orders 1024, 2048, and 4096 the relative factorization tolerance $20n\varepsilon$ is 0.625, 1.25, and 2.5 times the bfloat16 unit roundoff $2^{-8}$, and the second-round QR checker rejected only 49 of this control's 450 outputs.

\subsection{Per-entry QR results}
\label{sec:audit-qr-results}

Table~\ref{tab:qr-audit} reports both QR rounds.
Replay columns give evaluations failed under the first and second checkers, and discovery failures use the second checker throughout.
Under their own first-round checker, the five first-round entries
still failed 193, 73, 73, 245, and 3 of the 1,050 test cases in the
common grid.
The discovery counts combine the common grid with the threshold cases.
For example, the winner's 157 recorded failures comprise 155 from the common grid and two threshold cases.

\begin{table}[!htbp]
\centering
\small
\caption{QR results for the top five entries in each round and the controls.
Each target had 1{,}050--1{,}055 discovery evaluations and nineteen replay observations covering twelve distinct inputs.
QR discovery failures were not subjected to fresh-container confirmation.
The replay counts pool nineteen evaluations in which the seven dense inputs appear twice, so they differ from the counts over twelve inputs in Table~\ref{tab:qr-second-round}.
Data: audit record~\cite{melnichenko2026auditrecord}.}
\label{tab:qr-audit}
\renewcommand{\arraystretch}{1.15}
\begin{tabular}{@{}cclrrr@{}}
\toprule
Round & Rank & Entry & \makecell[br]{Replay failures,\\first checker} & \makecell[br]{Replay failures,\\second checker} & \makecell[br]{Discovery\\failures} \\
\midrule
1 & 1 & nikhilbarhate99 & 0/19 & 15/19 & 755 \\
1 & 2 & badelsteinlelbach & 0/19 & 0/19 & 92 \\
1 & 3 & Sinatras & 0/19 & 0/19 & 106 \\
1 & 4 & Maheshram1 & 3/19 & 15/19 & 586 \\
1 & 5 & elianaive & 0/19 & 0/19 & 3 \\
\midrule
2 & 1 & 10billiontokens & 0/19 & 0/19 & 157 \\
2 & 2 & gau.nernst & 0/19 & 0/19 & 46 \\
2 & 3 & dhu.randhar & 0/19 & 0/19 & 330 \\
2 & 4 & zhongmingee & 0/19 & 0/19 & 3 \\
2 & 5 & michaelmelons & 0/19 & 0/19 & 360 \\
\midrule
\multicolumn{3}{@{}l}{Baseline} & 0/19 & 0/19 & 0 \\
\multicolumn{3}{@{}l}{Baseline, bfloat16-rounded input} & 11/19 & 12/19 & 623 \\
\bottomrule
\end{tabular}

\end{table}

\subsection{Cholesky and eigendecomposition}
\label{sec:audit-other}

The audit ran 567 valid Cholesky inputs and 996 eigendecomposition inputs per entry.
A failure counts only when the organizers' checker rejects the entry's output and accepts the vendor routine's output on the same valid input.
Five failing candidates per entry were each tested in three reruns, two of them in fresh containers, and a candidate counts as reproduced only when all three fail and the vendor routine passes.
All five Cholesky entries and three of the five eigendecomposition entries had failures, and the failures of four and three entries were reproduced (Table~\ref{tab:audit-overview}).
One Cholesky entry failed 85 discovery inputs and passed all fifteen reruns, for a reason that remains unresolved.
The vendor eigensolver itself failed the checker on 38 inputs and raised a convergence error on two more, and these inputs are excluded from the entries' counts.
At large orders the eigendecomposition checker accepted outputs computed from input rounded to float8.
The audit record holds the per-entry counts~\cite{melnichenko2026auditrecord}.

\clearpage
\section{The library reuse analysis}
\label{app:benchmark_details}

\providecommand{\BenchTimeGain}{31.5\%}
\providecommand{\BenchTimeWins}{38}
\providecommand{\BenchCatOneTimeGain}{-2.6\%}
\providecommand{\BenchCatOneTimeLoss}{2.6\%}
\providecommand{\BenchCatTwoTimeGain}{56.8\%}
\providecommand{\BenchCatTwoTimeLoss}{56.8\%}
\providecommand{\BenchCatThreeTimeGain}{61.8\%}
\providecommand{\BenchCatThreeTimeLoss}{61.8\%}
\providecommand{\BenchZeroTimeGain}{4.3\%}
\providecommand{\BenchLastTimeGain}{65.7\%}
\providecommand{\BenchTokenGain}{22.1\%}
\providecommand{\BenchTokenWins}{34}
\providecommand{\BenchCatOneTokenGain}{-10.5\%}
\providecommand{\BenchCatOneTokenLoss}{10.5\%}
\providecommand{\BenchCatTwoTokenGain}{38.8\%}
\providecommand{\BenchCatTwoTokenLoss}{38.8\%}
\providecommand{\BenchCatThreeTokenGain}{56.7\%}
\providecommand{\BenchCatThreeTokenLoss}{56.7\%}
\providecommand{\BenchZeroTokenGain}{23.4\%}
\providecommand{\BenchLastTokenGain}{51.5\%}
\providecommand{\BenchCodeGain}{29.3\%}
\providecommand{\BenchCodeWins}{42}
\providecommand{\BenchCatOneCodeGain}{-1.8\%}
\providecommand{\BenchCatOneCodeLoss}{1.8\%}
\providecommand{\BenchCatTwoCodeGain}{41.7\%}
\providecommand{\BenchCatTwoCodeLoss}{41.7\%}
\providecommand{\BenchCatThreeCodeGain}{59.1\%}
\providecommand{\BenchCatThreeCodeLoss}{59.1\%}
\providecommand{\BenchZeroCodeGain}{3.9\%}
\providecommand{\BenchLastCodeGain}{54.2\%}
\providecommand{\BenchUsingTasks}{33}
\providecommand{\BenchZeroTasks}{27}

We evaluate the utility of \vscl{} in helping a coding agent formalize proofs from new research papers more efficiently. The evaluation focuses on results from Higham's \textit{Accuracy and Stability of Numerical Algorithms}~\cite{higham2002asna} that have been formalized in \vscl{}. Although the utility of the library might seem obvious at first glance, access to it does not necessarily improve efficiency. An agent not specifically trained on its contents may spend time and tokens searching for potentially useful results, assessing their relevance, and determining how to apply them. These costs could outweigh the savings from reusing existing library results instead of reconstructing them, particularly when no library result satisfying the three criteria described below is available to support a proof. We therefore examine whether existing library results can be reused in new, relevant proofs and whether access to the library reduces proof construction time, token consumption, and the amount of Lean code submitted. We also investigate how these benefits vary with expected library coverage and whether more extensive direct use of \vscl{} declarations is associated with greater reductions in these costs. Sections~\ref{subsec:benchmark_construct} to~\ref{sec:benchmark-results} describe the analysis and its results. The remaining subsections give the paper selection and task sources, the task categories, the task files and conditions, the agent protocol, the measurements, and the aggregation behind the figures. The code and scripts to reproduce this analysis are available in the GitHub repository ~\cite{qharabagh2026vsclreuse}.

\subsection{How is the analysis performed?}
\label{subsec:benchmark_construct}
We select 54 papers in numerical stability and related areas of numerical analysis which are relevant to the chapters of Higham's book. From these papers, we extract 60 proof formalization tasks, spanning topics such as floating-point arithmetic, matrix factorizations, iterative refinement, interpolation, and eigenvalue computations. We group the tasks according to the number of \vscl{} results judged applicable to their proofs that meet all three criteria: (i) reconstructing the result without access to \vscl{} would require significant mathematical effort, as judged by the LLM extracting the task (OpenAI Codex with GPT-6 Sol~\cite{openai2026codex}); (ii) the result goes beyond basic algebraic manipulations; and (iii) Mathlib does not already provide a readily usable equivalent. The categories are:
\begin{enumerate}
    \item \textbf{C1 (30 tasks):} No substantial library result is applicable.
    \item \textbf{C2 (20 tasks):} Exactly one substantial library result is applicable.
    \item \textbf{C3 (10 tasks):} At least two substantial library results are applicable.
\end{enumerate}

These criteria prevent tasks from being classified on the basis of elementary results, such as basic algebraic identities.
For example, \vscl ~contains an identity stating that subtraction distributes over finite expectations, $\mathbb{E}[X-Y]=\mathbb{E}[X]-\mathbb{E}[Y]$.
Although this identity may be useful as a minor step in a proof, its availability alone is not sufficient to place a task in C2 or C3. By contrast, \vscl ~contains a forward-error bound for floating-point matrix multiplication that bounds the error in each computed entry in terms of the magnitudes of the input matrix entries and a rounding-error factor. Reconstructing this bound requires a rounded-dot-product error analysis and supporting results, which the extracting LLM judges to involve significant mathematical effort. The bound also goes beyond basic algebraic manipulations, and Mathlib does not provide a readily usable equivalent. It therefore illustrates a library result satisfying all three criteria. Without this distinction, an extracting agent might classify a task as benefiting from library support merely because an elementary identity could appear in its proof, distorting the categories and weakening the practical significance of the subsequent measurements and conclusions.

For each task, we provide a \texttt{target.lean} file with the fixed  statement to be proved, a \texttt{context.md} file with the source statement and proof context, and a \texttt{Definitions.lean} module with definitions formalized from the paper. Each file is supplied unchanged and read-only in both conditions. In particular, the task statements do not depend on definitions from \vscl{}.

\subsection{Experimental setup}
\label{sec:benchmark-setup}
We compare a coding agent under two conditions: N (without \vscl{}) and L (with read-only access to \vscl{} and its precompiled Lean artifacts). Both receive the same task files and must prove the supplied statement without changing its assumptions, conclusion, or definitions. We use GPT-6 Sol with the \texttt{xhigh} reasoning setting, Lean v4.29.0-rc3, and Mathlib revision \texttt{e8ea1afc3279}. Both conditions have the same computational resources and may search local files, run shell commands, and compile Lean candidates. Neither has internet access or may invoke subagents.

Each of the 60 tasks is repeated three times in each condition, yielding 180 N/L pairs and 360 attempts, each with fresh agent state. Submissions must compile and pass integrity checks rejecting incomplete proofs, unauthorized axioms, and attempts to bypass Lean's proof checker.

We measure elapsed time from prompt release to submission, total input and output tokens (including cached input and reasoning tokens), and nonblank, noncomment lines of code (LOC) in the submitted \texttt{candidate.lean}. We also count \emph{direct proof declarations} (DP): distinct \vscl{} declarations referenced by the compiled target proof or reachable local helper bodies, excluding indirect dependencies within the library. 

For each task, we average resource costs over three repetitions. Figure~\ref{fig:benchmark-category-resources} compares category means between N and L; Figure~\ref{fig:benchmark-dp-gains} groups tasks by their mean DP in L and plots task-level gains. In each comparison, the signed percentage gain/loss is $g=100(1-X_L/X_N)$, where $X_N$ and $X_L$ are the relevant mean costs. Positive values indicate savings with L, and negative values indicate losses. For Figure~\ref{fig:benchmark-category-resources}, we calculate percentage gains from the mean resource costs within each category. For Figure~\ref{fig:benchmark-dp-gains}, we first calculate percentage gains for individual tasks and then average these gains within each DP group.

\subsection{Results}
\label{sec:benchmark-results}
All 360 attempts produced proofs that compiled successfully and passed the benchmark's integrity checks. Thus, both conditions achieved a 100\% success rate, and we focus on resource consumption.

\begin{figure}[p]
\centering
\includegraphics[width=.94\linewidth]{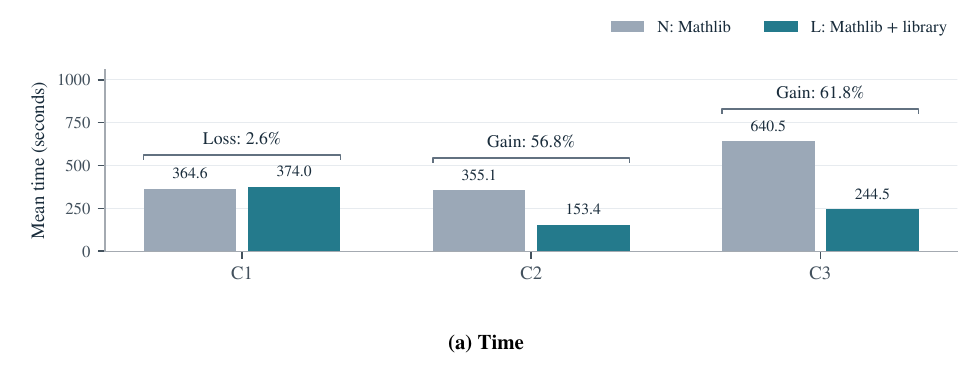}\par
\vspace{2pt}
\includegraphics[width=.94\linewidth]{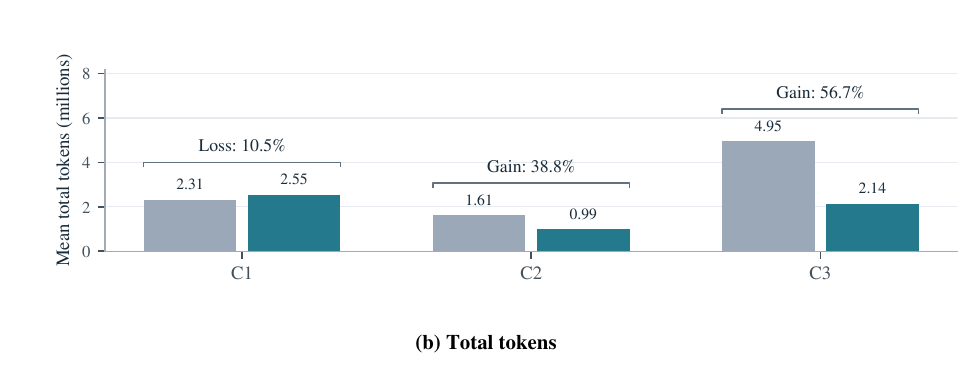}\par
\vspace{2pt}
\includegraphics[width=.94\linewidth]{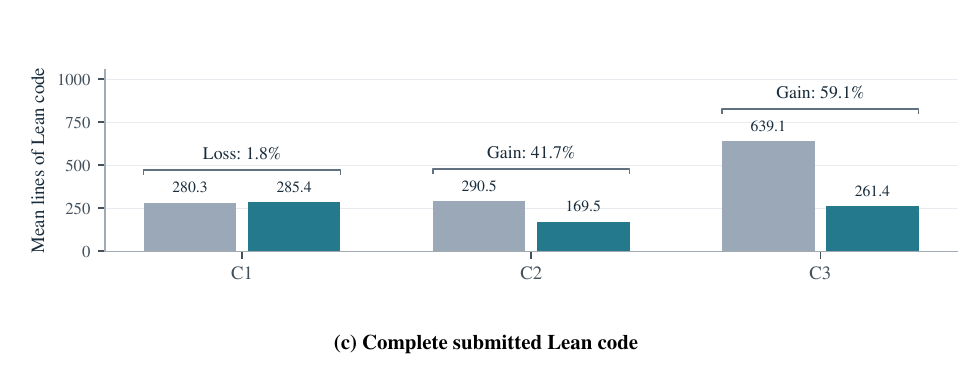}
\caption{ Resource consumption by expected library support.
Categories C1, C2, and C3 are defined by the number of results
available in \vscl{} that satisfy the three criteria and can be
used in constructing the task proofs: none, exactly one, and at least two, respectively.
For each category, bars show mean resource consumption per attempt
under conditions N (without \vscl{}) and L (with \vscl{}),
averaged across tasks and three repetitions per task.
Panels report (a) execution time, (b) total tokens, and
(c) lines of submitted Lean code.
Percentage gains are calculated as
$100(1-X_L/X_N)$, where $X_N$ and $X_L$ are the corresponding
category mean resource costs. Positive values indicate savings
with L; negative values are labelled as losses.
Total tokens include cached input, and code lines cover the
complete submitted Lean file.}
\label{fig:benchmark-category-resources}
\end{figure}

\paragraph{Results by expected library support.}
Figure~\ref{fig:benchmark-category-resources} shows that, for C2, L reduces time by \BenchCatTwoTimeGain, total tokens by \BenchCatTwoTokenGain, and code lines by \BenchCatTwoCodeGain. For C3, the corresponding reductions are \BenchCatThreeTimeGain, \BenchCatThreeTokenGain, and \BenchCatThreeCodeGain. In contrast, C1 shows increases in time, tokens, and code lines of \BenchCatOneTimeLoss, \BenchCatOneTokenLoss, and \BenchCatOneCodeLoss, respectively. The C1 overhead is consistent with the cost of searching the library when no result meeting the three criteria is available to offset that cost.

\begin{figure}[!htbp]
\centering
\includegraphics[width=.82\linewidth]{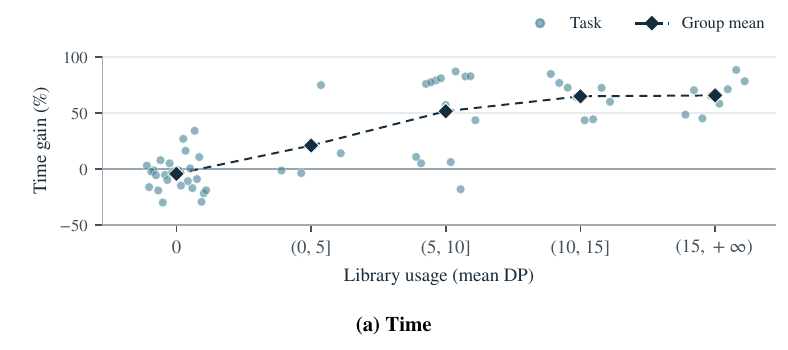}\par
\vspace{2pt}
\includegraphics[width=.82\linewidth]{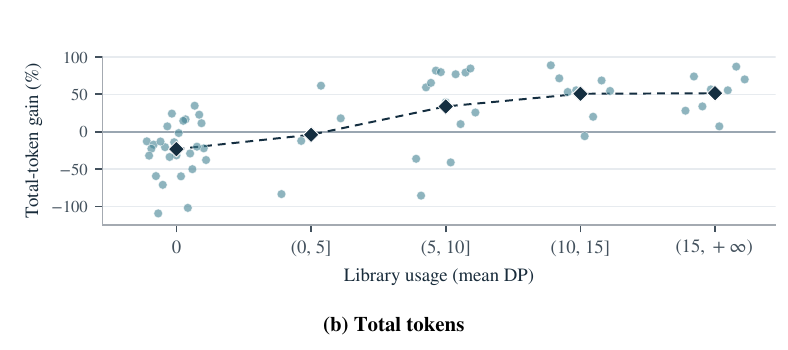}\par
\vspace{2pt}
\includegraphics[width=.82\linewidth]{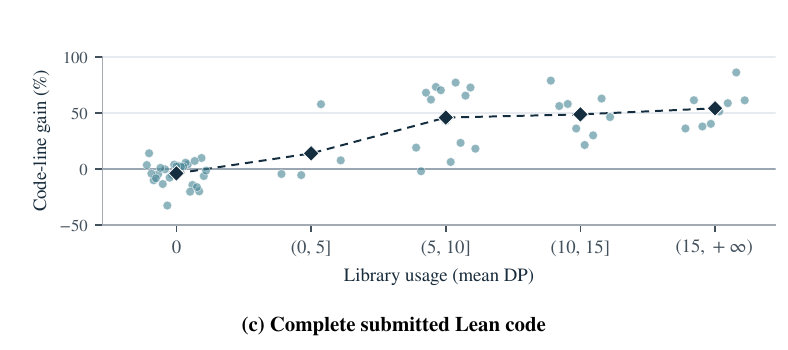}

\caption{ Gains in condition L relative to condition N, grouped by
library usage. DP (direct proof declarations) counts the distinct
library declarations used by the target proof and its local helper
proofs, excluding declarations used internally by imported library
proofs. For each task, DP is averaged over three L repetitions.
Tasks are grouped into intervals based on their mean DP, with zero
usage treated separately. Points represent individual task gains,
while diamonds indicate the mean gain across tasks within each group.
Dashed lines connect the group means as visual guides.
Positive gains indicate an advantage for L.}
\label{fig:benchmark-dp-gains}
\end{figure}

\paragraph{Gains with increasing use of \vscl{} declarations.}
The L proofs use library declarations on \BenchUsingTasks\ tasks; the remaining \BenchZeroTasks\ have zero DP in all repetitions. Mean DP is 0.34, 9.48, and 17.73 in C1, C2, and C3, respectively. Figure~\ref{fig:benchmark-dp-gains} groups tasks by observed DP. The zero-use group has mean losses of \BenchZeroTimeGain\ in time, \BenchZeroTokenGain\ in total tokens, and \BenchZeroCodeGain\ in lines of code. Across the displayed DP groups, mean gains generally increase with library use. For $\mathrm{DP}>15$, the mean gains reach \BenchLastTimeGain, \BenchLastTokenGain, and \BenchLastCodeGain\ in time, total tokens, and lines of code, respectively. For some lower-DP tasks, search overhead can outweigh the benefits of declarations that offer only limited help in constructing the proof, such as the elementary identity discussed in Section \ref{subsec:benchmark_construct}.

\paragraph{What the results suggest.}
The largest savings occur in C2 and C3, where library results satisfying the three criteria are expected to be reusable. L reduces proof time, token consumption, and lines of submitted Lean code, with greater gains in C3. Similarly, higher observed DP correlates with larger resource savings, supporting the effectiveness of \vscl{} when relevant results are available.

In C1, where no library results satisfying these criteria are expected to be available, N and L perform similarly, with L incurring slight search overhead. This highlights the need to expand \vscl{} coverage to reduce the proportion of C1 tasks and improve library search and access efficiency to minimize overhead.

\subsection{Paper selection and task sources}
We construct the task set by extracting proof formalization tasks from research papers. We use a task-driven selection process, screening a candidate collection of 78 papers in numerical stability and related areas of numerical analysis. A paper is eligible if it contains a result that can be formulated as a nontrivial, source-faithful Lean proof task, with sufficient statement and proof context to specify its mathematical assumptions and intended conclusion.

Our final selection comprises 54 papers published in leading venues in the field. These include 18 from \textit{SIAM Journal on Scientific Computing} and 18 from \textit{SIAM Journal on Matrix Analysis and Applications}. We also select 7 from \textit{IMA Journal of Numerical Analysis}, 5 from \textit{Mathematics of Computation}, and 3 from \textit{Numerische Mathematik}. We select one paper each from \textit{Advances in Computational Mathematics}, the IEEE Symposium on Computer Arithmetic, and the Vanderbilt approximation volume. The papers were selected so that their topics match chapters of Higham's book.

The tasks consist of either complete mathematical results or selected parts of results. They are drawn from theorems, lemmas, propositions, numbered equations, prose statements, or unnumbered mathematical displays. More specifically, 43 of the 60 tasks are drawn from explicitly labeled mathematical results. Of these, 21 are from theorems, 19 from lemmas, and 3 from propositions. The remaining 17 tasks are extracted from numbered equations, prose statements, or unnumbered mathematical displays.

The tasks span a broad range of topics. These include floating-point summation and dot products; mixed-precision matrix multiplication; triangular solves and matrix inversion; and LU, QR, and low-rank factorizations. They also include iterative refinement and Krylov methods; polynomial and rational interpolation; and eigenvalue computations. Further topics include matrix functions, Lyapunov equations, and probabilistic rounding-error analysis. 

\subsection{Task categories and expected coverage}
The classification is based on the number of substantial results available in \vscl{} that a coding agent could leverage when constructing a proof for the task. We consider a result in \vscl{} substantial if it satisfies three conditions: (i) reconstructing it without access to \vscl{} would require significant mathematical effort, as judged by the LLM extracting the task (OpenAI Codex with GPT-6 Sol~\cite{openai2026codex}); (ii) it goes beyond basic algebraic manipulations; and (iii) Mathlib does not already provide a readily usable equivalent. Thus, C1 contains 30 tasks with no applicable substantial result, C2 contains 20 tasks with exactly one, and C3 contains 10 tasks with at least two.

We use this categorization to cover different levels of support from \vscl{}. These range from tasks with no applicable substantial results to tasks that can use several library results. This allows us to examine how proof formalization performance varies with the library's coverage of the relevant mathematics.

The number of tasks in each category does not necessarily reflect the prevalence of that category in the numerical stability analysis literature. Nevertheless, we generally observe that, for a given collection of papers, C1 tasks are more common than tasks in the other two categories. Many C1 tasks, however, are not directly concerned with numerical stability. Instead, they consist of exact algebraic theorems, lemmas, and propositions. Authors use these as intermediate results to establish other claims. C2 is generally the second most common category. C3 tasks occur less frequently than C2 tasks.

\subsection{Task files and experimental conditions}
For each task, we create a \texttt{target.lean} file containing the formal statement to be proved and a \texttt{context.md} file providing its natural-language context, including the informal statement, proof, and relevant information from the source paper. We also create a \texttt{Definitions.lean} module containing formalized definitions from the corresponding paper. These definitions are used to express the formal statement of the task.

To evaluate library reuse, we consider two experimental conditions: one in which the coding agent has access to \vscl{} and another in which it does not. In both conditions, the \texttt{target.lean}, \texttt{context.md}, and \texttt{Definitions.lean} files are provided to the agent as fixed, read-only files that cannot be modified. Providing the same \texttt{target.lean} file in both conditions ensures that the agents work on identical formalization tasks, enabling a fair comparison. The \texttt{Definitions.lean} module makes this possible by ensuring that the formal statements do not depend on definitions from \vscl{}.

We denote the conditions N (without \vscl{}) and L (with \vscl{}). The agent must prove the supplied statement without modifying its assumptions, conclusion, or supporting definitions. Both conditions allow shell commands and interaction with the Lean environment. In addition, L has read-only access to \vscl{}, mounted together with its precompiled Lean artifacts (\texttt{.olean} files), so library modules can be imported without rebuilding them from source.

\subsection{Model, resources, and agent protocol}
We use GPT-6 Sol with the \texttt{xhigh} reasoning setting. Both conditions use Lean v4.29.0-rc3 and the same Mathlib revision, \texttt{e8ea1afc3279}.
Both conditions run on the same hardware configuration, with four physical cores of an Intel Core i9-13900K processor and a 16 GiB memory limit per agent attempt. CPU and memory allocations are enforced to ensure consistent computational resources across conditions.

In both conditions, the agent receives a prompt instructing it to construct a proof of the fixed statement in \texttt{target.lean} and submit its final code in \texttt{candidate.lean}. The prompt in condition L contains an additional section specifying the path to \vscl{}, which contains a fixed snapshot of the material formalized from \textit{Accuracy and Stability of Numerical Algorithms}~\cite{higham2002asna}, and informing the agent that it may use results from the library when needed. Agents may search their mounted files using tools such as \texttt{rg} and \texttt{find}, inspect local library sources, and compile Lean candidates while developing their proofs. Agents are not allowed to invoke subagents and have no internet access. Internet access is disabled to prevent agents from retrieving external solutions, including online copies of \vscl{}, while subagents are prohibited to prevent differences in parallelism from affecting resource comparisons between conditions.

Each of the 60 tasks discussed in Section~\ref{subsec:benchmark_construct} is repeated three times in each condition, yielding 180 N/L pairs and 360 attempts. Every attempt starts with a fresh workspace, conversation, and local agent state. A submitted proof must compile successfully and pass the integrity checks, which reject incomplete proofs, unauthorized axioms, and attempts to bypass Lean's proof checker.

\subsection{Resource measurements and direct library use}
We measure execution time from the release of the prompt to the final submission, including local searches and interactions with Lean. Post-submission validation is outside this interval. Token consumption is the sum of all input tokens, including cached input, and output tokens. Cached input is part of the input total, and reasoning tokens are part of the output total; neither is added a second time. We also measure the lines of code (LOC) in the agent's final submitted solution, \texttt{candidate.lean}. This count includes nonblank, noncomment code lines, including import statements and local helper declarations. Lines belonging to declarations in imported library modules are excluded. LOC serves as a proxy for the amount of proof and helper code submitted locally. Lower LOC can reflect library reuse, shorter tactics, or a different proof strategy; it does not by itself establish reuse.

As a more direct measure of library use, we also count \emph{direct proof declarations} (DP). These are the distinct \vscl{} declarations referenced by the compiled target proof or by the bodies of local helpers reachable from that proof in the submitted file. A declaration referenced multiple times is counted only once. If a library theorem internally uses another library result, that result is not counted unless it is also directly referenced by the submitted proof or a reachable local helper. The analysis uses resolved Lean declarations, allowing it to detect references even when namespace prefixes are omitted. The count includes definitions and other constants as well as theorems.

\subsection{Aggregation and figure interpretation}
For each task and condition, we average resource consumption over the three repetitions. Figure~\ref{fig:benchmark-category-resources} compares agent performance under conditions N and L across the coverage-based categories C1--C3. Each bar represents the mean resource consumption across the tasks in the corresponding category. For category mean resource costs $X_N>0$ and $X_L\geq0$ in conditions N and L, respectively, the signed percentage gain/loss is
\[
g=100\left(1-\frac{X_L}{X_N}\right).
\]
A positive value indicates a gain (lower resource use with L), while a negative value indicates a loss (higher resource use with L). The bar brackets label losses using their positive magnitudes. For example, $g=-2.6$ is shown as ``Loss: 2.6\%.''

Figure~\ref{fig:benchmark-dp-gains} shows the percentage gains and losses under condition L relative to condition N as a function of direct library use, measured by the DP count. We average the DP counts across each task's three repetitions in condition L and assign the task to one of five groups based on this mean count: zero, $(0,5]$, $(5,10]$, $(10,15]$, or $(15,+\infty)$. Each point represents one task, with its signed percentage gain/loss calculated as
\[
g_{\mathrm{task}}=100\left(1-\frac{\bar{x}_L}{\bar{x}_N}\right),
\]
where $\bar{x}_N$ and $\bar{x}_L$ denote the mean resource consumption across the three repetitions of the task in conditions N and L, respectively. Each denominator is positive for the reported measurements. Diamonds indicate the equally weighted mean of these signed task-level values within each group, while dashed lines connect the means as a visual guide. Positive values indicate gains and negative values indicate losses; the scatter plots retain these signs.

The DP counts are not a direct measure of the mathematical significance of reused results: some counted declarations may be general-purpose helpers that do not substantially reduce proof construction effort. Therefore, the overhead of accessing and searching the library can outweigh the benefits for some tasks in lower-DP groups.

\section{Open-weight models: costs, experiments, and workflow}
\label{sec:open-weights}

Formalizing material from books requires months of sustained
model use, as each claim may require several rounds of source
interpretation, statement review, proof generation, and repair.
Unsuccessful attempts also consume resources. Extending \vscl{} to a
broad collection of books therefore requires a substantial investment
in model access and human supervision.

A simple estimate illustrates the expenses. Based on our
experience, formalizing selected results from each advanced mathematics
book typically requires approximately two to three
months using an agent such as Codex, assuming that human intervention
is limited to cases in which the agent flags uncertainty about how to
proceed. Under this assumption, six subscriptions at US\$200 per month
for two to three months would cost US\$2,400--US\$3,600 in subscription fees
alone~\cite{openai2026codexpricing}.
This is just an illustrative budget for library development. It excludes
human supervision, mathematical review, and other computing costs.
Personal subscriptions also impose usage limits that constrain
sustained execution. API access supports larger workloads, but charges
accumulate across the full sequence of model
calls~\cite{openai2026apipricing} and long contexts, independent reviews,
and repeated proof attempts all contribute to this cost. An observed
API run provides a reference for the potential expense. One of our 13-hour runs cost approximately US\$1,300, corresponding
to an average of US\$100 per hour. Sustaining this spending rate continuously for
thirty days would cost approximately US\$70,000 in
API charges alone. The actual expense depends on the
model settings, number of concurrent agents, token consumption, and
caching. Human supervision, mathematical review, and other computing
costs are additional. These expenses motivate methods that reduce the
cost of extending \vscl{} while preserving its formalization capability.

Open-weight models offer one possible approach. Their weights can be
deployed on shared computing infrastructure, where we control
scheduling and resource allocation. We can also adapt the models to
recurring tasks and use specialized models for proof generation and
repair. Models such as Leanstral~\cite{mistral2026leanstral} provide publicly
available weights for proof engineering in Lean. With
sufficient use of the hardware and adequate formalization performance,
this approach may reduce the cost of processing a large collection
of claims.

The relevant comparison is the total cost per accepted formalization.
It must include computing, maintenance, training, and review. Lower
inference costs may be offset by additional proof attempts or greater
human involvement. Our goal is therefore to reduce the cost of library
construction while maintaining source faithfulness, proof completion,
and coverage.

\subsection{Our current approach to open-weight models}

We are pursuing two directions. The first is to run the formalization
skill described in Section~\ref{sec:skill} using a self-hosted,
open-weight coding agent instead of OpenAI's Codex~\cite{openai2026codex}.
The second is to develop an orchestrator--prover workflow that
separates statement construction and review from proof generation.

\paragraph{Standalone formalization experiments.}
To investigate the use of open-weight models for extending \vscl{}, we ran separate
formalization experiments with DeepSeek-V4.1-Flash-1017~\cite{deepseekai2026deepseekv4}
through OpenRouter and with Codex, configured with
GPT-6.1 Sol at the high reasoning setting. Each model attempted the four-page
Appetizer and seven-page Chapter~1 of \emph{High-Dimensional
Probability}~\cite{vershynin2018hdp}. Each experiment began with the source PDF
and a fresh Lean project. The agents could use Lean and the same fixed version
of Mathlib, but could not access \vscl{} or earlier formalizations of the material.
They therefore had to construct the formal statements as well as their proofs.

DeepSeek-V4.1-Flash-1017 followed a simpler source-to-Lean protocol;
its Chapter~1 run included four operator-requested continuations. GPT-6.1 Sol used an
adapted version of the formalization skill, with separate agents reviewing
statements and checking proofs. Both GPT-6.1 Sol runs stopped to request clarification
of the source. For the Appetizer, we supplied a precise lower-bound statement
for the Notes' assertion that the covering bound is sharp, based on an entropy
estimate of Carl and Pajor~\cite{Carl1988GelfandNO}.
For Chapter~1, we specified positive integer orders for raw moments of signed
random variables and positive real orders for raw moments of nonnegative
variables, while retaining positive real orders for absolute moments. GPT-6.1 Sol
then resumed the existing work. Its final results are therefore assisted
formalizations of this amended scope.

Both developments contained substantial proofs of the main approximation and
probability limit results. The Lean projects built successfully, and their
axiom checks passed. Separate model reviews compared the formal statements and
their definitions with the PDFs. DeepSeek-V4.1-Flash-1017 had not proved the Appetizer's final
covering-sharpness claim. Some of its expectation and moment definitions also
required domain qualifications because their real-valued integrals do not
represent infinite or undefined expectations. GPT-6.1 Sol provided more explicit
treatment of these cases and proved the supplied sharpness statement. Its audit
found no blocking correspondence error within the recorded amended scope.
Both developments nevertheless retain finite-second-moment assumptions for
variance results. Neither audit establishes an unrestricted, literal
formalization of every source statement. As in Section~\ref{sec:skill}, these source
comparisons are judgments by reviewing agents, separate from Lean's verification
of the proofs.

Table~\ref{tab:standalone-hdp} in Section~\ref{sec:open-weights-summary}
reports the observed resource use. Active time
includes model-side setup, local searches, Lean commands, and continuations,
but excludes waiting between turns for operator input. Operator preparation,
specification research, and the final audits are outside these measurements.
Reasoning tokens are included in output, and the GPT-6.1 Sol totals include the
assisted continuation for each unit.
The DeepSeek-V4.1-Flash-1017 costs are recorded API charges. GPT-6.1 Sol ran through ChatGPT sign-in~\cite{openai2026codexpricing},
so its API costs are estimates calculated from recorded tokens, observed
caching, and the applicable API rates~\cite{openai2026apipricing}, not charges allocated from the
subscription.

These experiments show that a hosted open-weight model can construct substantial
Lean proofs from textbook material. They also show why
proof validity, source correspondence, and coverage must be assessed separately.
The differing procedures, assistance, and final scopes prevent a controlled
comparison of model capability or cost per accepted formalization. These results
support further investigation of open-weight models for extending the library.

\begin{figure}[h!]
\centering
\includegraphics[width=\linewidth]{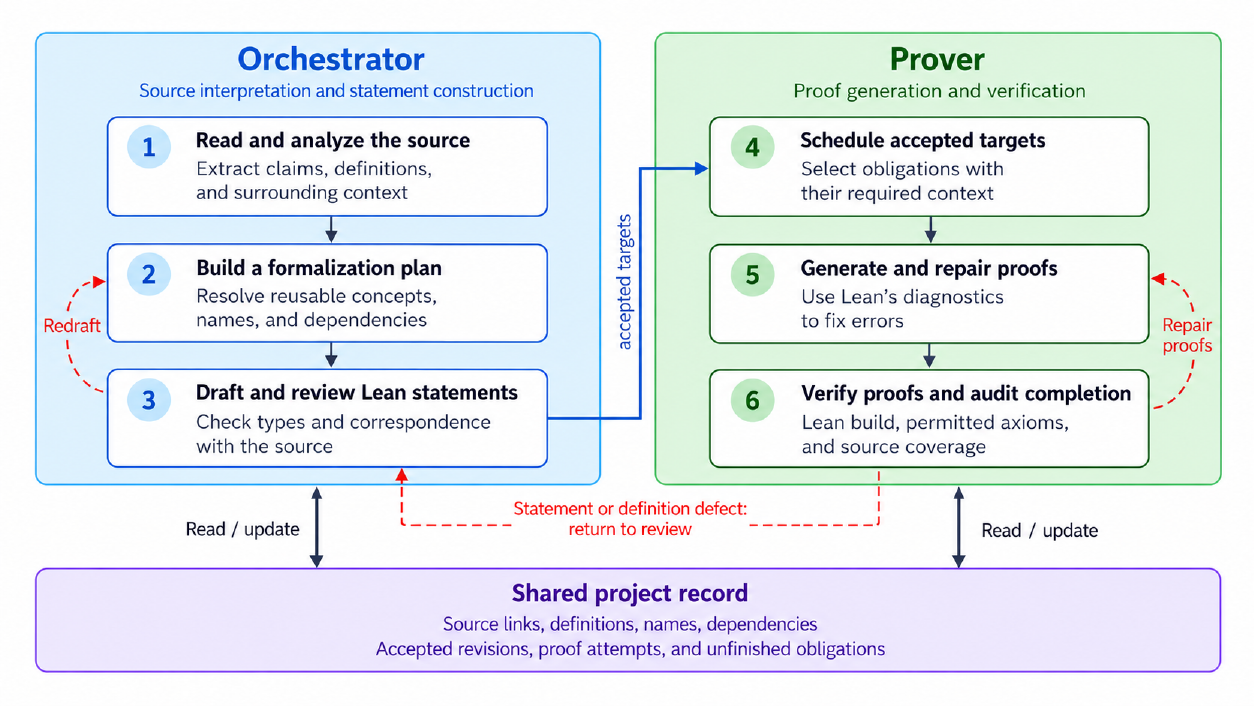}
\caption{The proposed orchestrator--prover workflow separates statement construction
and review from proof generation. The orchestrator prepares Lean statements
and arranges their review against the source before passing accepted
targets to the prover. Proof errors return to repair, while errors in
statements or definitions return to review. Both roles use a shared
project record. Acceptance requires Lean verification and software checks.
The full workflow remains under development.}
\label{fig:open-weight-phases}
\end{figure}

\paragraph{The orchestrator--prover workflow.}
We aim to reduce the cost of self-hosted autoformalization by using two
smaller, specialized open-weight models within a predefined workflow,
rather than a coding agent powered by a large model that organizes
and carries out the entire process using only skill files as guidance.
Running such a large model locally can require substantial computing
resources. By defining the roles, stages, and checks in advance, we
aim to enable smaller models to perform their assigned tasks without
having to manage the entire process autonomously.
The \emph{orchestrator} prepares Lean definitions and theorem statements
from the source text and coordinates their review. The \emph{prover}
constructs proofs of the reviewed statements
(Figure~\ref{fig:open-weight-phases}).

In steps~1--3, the orchestrator identifies each claim's assumptions,
definitions, and supporting passages. It checks for reusable library
definitions and results before drafting new declarations, and records
the dependencies between claims. Lean checks that the definitions and
statements are well typed; a separate review checks that they match the
source. For example, the covering result discussed above must require
every point of the set to lie in at least one of the chosen balls.
Merely introducing the balls would not express the source result, even
if Lean accepted the statement as well typed.

Software schedules the reviewed statements for proof (step~4).
The prover then writes proofs and uses Lean's error messages to repair
unsuccessful attempts (step~5). It may change the proof code but not
the reviewed statements or definitions. If a statement or definition
is incorrect, the claim returns to the orchestrator for correction
and another review. Finally, Lean verifies the proofs, while audit
software checks that they use only permitted axioms and records which
source claims have been proved and which remain unfinished (step~6).

Both roles use a shared project record containing source links,
accepted definitions and statements, dependencies, proof attempts,
and unfinished work. This allows a new worker to continue without
needing the previous conversation. The separation of roles and use
of a shared record follow the organization shown in
LeanFlow~\cite{milikic2026leanflow}.

\paragraph{Development of the orchestrator--prover workflow.}
The source-processing and planning software for this workflow has been implemented and tested
on Vershynin's \emph{High-Dimensional Probability}~\cite{vershynin2018hdp}.
It extracts and indexes the text, records dependencies and links planned
declarations to the passages they are meant to formalize. These records
support the later construction of Lean declarations, but no real formalization work has
been done so far. Statement construction, review, proof repair, and final acceptance still
need to be integrated and tested together.

\paragraph{Evaluation of the orchestrator--prover workflow.}
The next step for the orchestrator--prover workflow is to obtain a small result with a reviewed statement and a
complete Lean proof before testing larger held-out chapters. These tests will
use matched environments and budgets, with translation and proof repair
assessed separately. They will examine whether shared planning and the
checks between stages reduce repeated failures and human intervention.
Where the base models fall short, we will investigate adaptation using
reviewed translations and verified repairs.

The evaluation will measure source faithfulness and the proportion of claims
fully proved, alongside runtime, human intervention, and compute cost.
Training costs will be reported separately. Lean acceptance establishes proof
validity, while the judgment that a statement expresses the source claim
remains a separate review. Both open-weight model directions remain under investigation.

\end{document}